\documentclass[12pt]{article}
   \usepackage{style-perso}
\usepackage{amsmath}
\usepackage{amsthm}
\usepackage{amssymb}
\usepackage{amsfonts,amssymb,verbatim}
\usepackage{epsfig}
\usepackage{pb-diagram}

\input colordvi
\newtheorem{theorem}{Theorem}[section]

\newtheorem{lemma}{Lemma}[section]
\newtheorem{proposition}[theorem]{Proposition}
\newtheorem{definition}{Definition}[section]

\def\ovr{\overrightarrow}
\def\cT{\mathcal{T}}
\def\cA{\mathcal{A}}
\def\cC{\mathcal{C}}
\def\cD{\mathcal{D}}
\def\wt{\widetilde}
\def\cB{\mathcal{B}}
\def\cE{\mathcal{E}}
\def\cS{\mathcal{S}}

\def\cM{\mathcal{M}}
\def\ul{\underline}
\def\cP{\mathcal{P}}
\def\cH{\mathcal{H}}
\def\cR{\mathcal{R}}

\def\MVB{M'_{\ \!B}}
\def\MVBf{M'_{\ \!Bf}}
\def\MVf{M'_{\ \!f}}
\def\cG{\mathcal{G}}
\def\cZ{\mathcal{Z}}
\def\tC{\widetilde{C}}
\def\Seq{\mathrm{Seq}}
\def\Set{\mathrm{Set}}
\def\Cyc{\mathrm{Cyc}}
\def\cBv{\cB_{\circ}}
\def\cBe{\cB_{\circ-\circ}}
\def\cBh{\cB_{\circ\rightarrow\circ}}
\def\cBthree{\cB}
\def\cBR{\cB_R}
\def\cBM{\cB_M}
\def\cBT{\cB_T}
\def\cBRM{\cB_{R-M}}
\def\cBRT{\cB_{R-T}}
\def\cBMT{\cB_{M-T}}
\def\cBTT{\cB_{T-T}}
\def\cBRtoM{\cB_{R\rightarrow M}}
\def\cBMtoR{\cB_{M\rightarrow R}}
\def\cBRtoT{\cB_{R\rightarrow T}}
\def\cBTtoR{\cB_{T\rightarrow R}}
\def\cBMtoT{\cB_{M\rightarrow T}}
\def\cBTtoM{\cB_{T\rightarrow M}}
\def\cBTtoT{\cB_{T\rightarrow T}}
\def\cGcp{\cG_1\ \!\!\!'}
\def\cGbp{\cG_2\ \!\!\!'}
\def\cGbr{\ovr{\cG_2}}
\def\cGtp{\cG_3\ \!\!\!'}

\def\Mcr{\ovr{M}}
\def\Mbr{\ovr{L}}
\def\Mtr{\ovr{K}}

\def\Mcp{M'}
\def\Mbp{L'}
\def\Mtp{K'}

\def\uV{V}

\def\cV{\mathcal{V}}

\def\uD{{D}}
\def\uS{{S}}
\def\uP{{P}}
\def\uH{{H}}

\def\loga{\mathrm{loga}}

\def\ds{\displaystyle}

\def\ovr{\overrightarrow}

\title{A Complete Grammar for Decomposing a Family of Graphs into 3-connected Components}
\author{Guillaume Chapuy$^{1}$, \'Eric Fusy$^{2}$, Mihyun Kang$^{3}$ and Bilyana Shoilekova$^{4}$}

\begin{document}
\maketitle
\begin{abstract}
Tutte has described in the book ``Connectivity in graphs'' a canonical decomposition of any graph
into 3-connected components. In this article we translate (using the language of symbolic combinatorics)
 Tutte's decomposition into 
a \emph{general grammar} expressing any family $\cG$ of graphs (with some stability conditions)
in terms of the subfamily $\cG_3$ of graphs in $\cG$ that are 3-connected (until now, such a general grammar was only known for the decomposition into $2$-connected components).
As a byproduct, our grammar yields an explicit system of equations to express 
the series counting a (labelled) family of graphs in terms of the series counting the subfamily of 
$3$-connected graphs. A key ingredient we use is an extension of the so-called dissymmetry theorem, 
which yields negative signs in the grammar and associated equation system, but  
has the considerable advantage of avoiding the difficult \emph{integration steps} that appear 
with other approaches, in particular in recent work by Gim\'enez and Noy on counting planar graphs.

As a main application we recover \emph{in a purely combinatorial way} the analytic expression found
by Gim\'enez and Noy for the series counting labelled planar graphs (such an expression is crucial to do asymptotic enumeration  and to obtain limit laws of various parameters on random planar graphs).
Besides the grammar, 
an important ingredient of our method is a recent bijective construction of planar maps by Bouttier,
Di Francesco and Guitter.

Finally, our grammar applies also to the case of \emph{unlabelled} structures, since the dissymetry theorem takes symmetries into account. 
Even if there are still difficulties in counting unlabelled 3-connected planar graphs, 
we think that our grammar is a promising tool toward the asymptotic enumeration of unlabelled planar graphs,
since it circumvents some difficult integral calculations.
\end{abstract}

\noindent\rule{10cm}{1pt}

{\footnotesize
\noindent $^{1}$: LIX, \'Ecole Polytechnique, Paris, France. \texttt{chapuy@lix.polytechnique.fr}\\
$^{2}$: Dept. Mathematics, UBC, Vancouver, Canada. \texttt{fusy@lix.polytechnique.fr}\\
$^{3}$: Institut f\"ur Informatik, Humboldt-Universit\"at zu Berlin, Germany. \texttt{kang@math.tu-berlin.de}\\
$^{4}$: Department of Statistics, University of Oxford, UK. \texttt{shoileko@stats.ox.ac.uk}
}

\section{Introduction}
Planar  graphs and related families of structures have recently received a lot of attention both from a probabilistic and an enumerative point of view \cite{BeGa,BKLM07,Fu05a,gimeneznoy,MSW05}. 
While the probabilistic approach already yields significant qualitative results, the enumerative approach provides a complete solution
regarding the asymptotic behaviour of many parameters on random planar graphs (limit law for the number of edges, connected components), as demonstrated by Gim\'enez 
and Noy for planar graphs~\cite{gimeneznoy} building on earlier work of Bender, Gao, Wormald~\cite{BeGa}. Subfamilies of labelled planar graphs have been treated 
in a similar way in~\cite{BoGiKaNo07,BKLM07}.

The main lines of the enumerative method date back to Tutte \cite{Tu63,Tutte}, where graphs are decomposed into components of higher connectivity: A graph is decomposed into connected components, each of which is decomposed into 2-connected components, each of which is further decomposed into 3-connected components. For planar graphs every 3-connected graph has a unique embedding on the sphere, a result due to Whitney~\cite{Whi32}, hence the number of 3-connected planar \emph{graphs} can be derived from the number of 3-connected planar \emph{maps}.
This already makes it possible to get a polynomial time method for exact counting (via recurrences that are derived for the counting coefficients) 
and uniform random sampling of labelled planar graphs, as described
by Bodirsky et al~\cite{bodirsky}. 
This decomposition scheme can also be exploited to get asymptotic results: asymptotic enumeration, limit laws for various parameters.
In that case, the study is more technical and relies on two main steps: symbolic and analytic.
In the symbolic step, Tutte's decomposition is translated into an equation system satisfied by the counting series.
In the analytic step, 
a careful analysis of the equation system makes it possible to locate and determine the nature of the (dominant) singularities  of the counting series; from there,
transfer theorems of singularity analysis, as presented in the forthcoming book by Flajolet and Sedgewick~\cite{FlaSe}, yield the asymptotic results.

In this article we focus on the symbolic step: how to translate Tutte's decomposition into an equation system in an automatic way.
Our goal is to use a formalism as general as possible, which works both in the labelled and in the unlabelled framework, and works
 for a generic family of graphs (however under a certain stability condition), not only planar graphs.
Our output is a generic decomposition grammar---the grammar
 is shown in Figure~\ref{GrammarFig}---that corresponds to the translation of Tutte's decomposition. 
Getting such a grammar is however nontrivial, as Tutte's decomposition is rather involved; we exploit the \emph{dissymmetry theorem}
(Theorem~\ref{DissymTreesDecomp}) 
applied to trees 
that are naturally associated with the decomposition of a graph. Similar ideas were recently independently described by Gagarin et al in~\cite{GLLW08}, 
where they express a species of 2-connected graphs in terms of the 3-connected subspecies.
Translating the decomposition into a grammar as we do here is very transparent and makes it possible to easily get equation systems in an automatic way,
both in the \emph{labelled} case (with generating functions) and in the \emph{unlabelled} case (with P\'olya cycle index sums).
Let us also mention that, when performing the symbolic step in~\cite{gimeneznoy}, 
Gim\'enez and Noy also translate Tutte's decomposition into a positive equation system,  
but they do it only \emph{partially}, as some of the generating functions in the system they obtain have to be \emph{integrated}; therefore 
they have to deal with complicated analytic integrations, see~\cite{gimeneznoy} and more recently~\cite{GNR07} for a generalized presentation. 
In contrast, in the equation system derived from our grammar, no integration step is needed;
and as expected, the only terminal series are those counting the 3-connected subfamilies (indeed, 3-connected graphs are the terminal bricks
in Tutte's decomposition). In some way, the dissymmetry theorem used to write down the grammar allows us to do the integrations \emph{combinatorially}.

In addition to the grammar, an important outcome of this paper is to show that the analytic (implicit) expression for the series counting labelled planar graphs
can be found in a completely combinatorial way (using also some standard algebraic manipulations), thus providing an alternative more direct way compared to the method
of Gim\'enez and Noy, which requires integration steps. Thanks to our grammar, finding  an analytic expression for the series counting planar graphs 
reduces to finding one for the series
counting 3-connected planar graphs, which is equivalent to the series counting 3-connected maps by Whitney's theorem. 
Some difficulty occurs here, as only an expression for the series counting \emph{rooted} 3-connected maps
is accessible in a direct combinatorial way. So it seems that some integration step is needed here, and actually that integration was analytically solved by 
Gim\'enez and Noy in~\cite{gimeneznoy}. In contrast we aim at finding an expression for the series counting \emph{unrooted} 3-connected maps in a more direct 
combinatorial way. We show that it is possible, 
by starting from a bijective construction of vertex-pointed maps---due to Bouttier,
 Di Francesco, and Guitter~\cite{BoDiGu04}---and going down
to vertex-pointed 3-connected maps; then Euler's relation makes it possible to obtain the series counting 3-connected maps from the series
counting vertex-pointed and rooted ones. In some way, Euler's relation can be seen as a generalization of the dissymmetry theorem that applies to maps
and allows us to integrate ``combinatorially'' a series of rooted maps.

Concerning unlabelled enumeration, we prefer to stay very brief in this article (the counting tools are cycle index sums, which are 
a convenient refinement of ordinary generating functions). Let us just mention that our grammar can be translated into a generic equation system
relating the cycle index sum (more precisely, a certain refinement w.r.t. edges) of a family of graphs to the cycle index sum of the 3-connected subfamily. However such a system is very complicated. Indeed the 
relation between 3-connected and 2-connected graphs involves edge-substitutions, which are easily addressed by exponential generating functions for labelled enumeration  (just 
substitute the variable counting edges) but are more intricate when it comes to unlabelled enumeration (the computation rule is a specific multivariate substitution). 
We refer the reader to the recent articles by Gagarin et al~\cite{GLL07,GLLW08}  for more details.
And we plan to investigate the unlabelled case in future work, 
in particular to recover (and possibly extend) in a unified framework  
 the few available results on counting asymptotically unlabelled subfamilies of planar graphs~\cite{Sh07,BoFuKaVi07b}. 

\vspace{.2cm}

\textbf{Outline.} After the introduction, 
there are four preliminary sections to recall important results in view of writing down the grammar. 
Firstly we recall in Section \ref{sec:Symbolic}  the principles of the symbolic method, which makes it possible to translate systematically combinatorial decompositions into enumeration results, using generating functions for labelled classes and ordinary generating functions (via cycle index sums) for unlabelled classes. 
In Section \ref{DissymThms} we recall the dissymmetry theorem for trees and state an extension of the theorem to so-called tree-decomposable classes. In Section \ref{DeGr3-cComps} we give an outline of the necessary graph theoretic concepts for the decomposition strategy. Then we recall the decomposition of connected graphs into 2-connected components and of 2-connected graphs into 3-connected components, following the description of Tutte~\cite{Tutte}. We additionally give precise characterizations of the different trees resulting from the decompositions. 

In the last three sections, we present our new results.
In Section \ref{sec:Decomp Grammar} we write down the grammar resulting  from 
Tutte's decomposition, thereby making an extensive use of the dissymmetry theorem. The complete 
grammar is shown in Figure~\ref{GrammarFig}.
In Section~\ref{sec:applica}, we discuss applications to labelled enumeration; the grammar is translated into an equation system---shown in Figure~\ref{fig:system}---expressing a series counting a 
graph family in terms of the series counting the 3-connected subfamily. 
 Finally, building on this and on enumeration techniques for maps, 
 we explain in Section~\ref{sec:planar}  how to get an (implicit)  
 analytic expression for the series counting labelled planar graphs.

\section{The symbolic method of enumeration}\label{sec:Symbolic}

In this section we recall important concepts and results in symbolic combinatorics, 
which are 
presented in details in the book by Flajolet and Sedgewick~\cite{FlaSe} (with an emphasis
on analytic methods and asymptotic enumeration) and the book 
by Bergeron, Labelle, and Leroux~\cite{BeLaLe} (with an emphasis on unlabelled enumeration). 
The symbolic method is a theory for enumerating \emph{decomposable} combinatorial classes in a systematic way. The idea is to find a recursive decomposition for a class $\cC$, and to write this decomposition as a \emph{grammar} involving a collection of \emph{basic classes} and \emph{combinatorial constructions}. The grammar in turn translates to a recursive equation-system 
satisfied by the associated generating function $C(x)$, which is a formal series whose coefficients are formed from the counting sequence of the class $\cC$.
From there, the counting coefficients of 
$\cC$ can be extracted, either in the form of an estimate (asymptotic enumeration),
or in the form of a counting process (exact enumeration).

\subsection{Labelled/unlabelled structures}
 A \emph{combinatorial class} $\cC$ (also called a species of combinatorial structures) 
 is a set of \emph{labelled} objects equipped with a size function; each object of $\cC$ is made of $n$ atoms (typically, vertices of graphs) assembled in a specific way, the atoms bearing distinct labels in $[1..n]:=\{1,\ldots, n\}$ (in the general theory of species, any system of labels is allowed). The number of objects of each size $n$, denoted $\cC_n$, is finite. The  classes we consider are stable under isomorphism (two structures  are called \emph{isomorphic} if one is obtained from the other by relabelling the atoms). 
 Therefore, the labels on the atoms only serve to distinguish them, which means that no notion of order is used for the labels. The class of objects in $\cC$ taken up to 
  isomorphism is called the \emph{unlabelled class} of $\cC$ and
is  denoted by $\tC=\cup_n\tC_n$. 

\subsection{Basic classes and combinatorial constructions}
\label{subsec:basic}
We introduce the basic classes and combinatorial constructions, as well as the rules to compute the associated counting series. The neutral class $\cE$ is made of a single object of size $0$. The atomic class $\mathcal{Z}$ is made of a single object of size $1$.  Further basic classes are the $\Seq$-class, the $\Set$-class, and the $\Cyc$-class, each object of the class being a collection of $n$ atoms assembled respectively as an ordered sequence, an unordered set, and an oriented cycle.

Next we turn to the main constructions of the symbolic method. The \emph{sum} $\cA+\cB$ of two classes $\cA$ and $\cB$ refers to the disjoint union of the classes.
The \emph{partitional product} (shortly product) $\cA*\cB$ of two classes $\cA$ and $\cB$ is the set of labelled objects that are obtained as follows: take a pair $(\gamma\in\cA,\beta \in\cB)$, distribute distinct labels on the overall atom-set (i.e., if $\beta$ and $\gamma$  are of respective sizes $n_1, n_2$, then the set
of labels that are distributed is $[1..(n_1+n_2)]$), and forget the original labels on $\beta$ and $\gamma$. 
Given two classes $\cA$ and $\cB$ with no object of size $0$ in $\cB$, the \emph{composition}  of $\cA$ and $\cB$, is the class $\cA\circ\cB$ ---also written $\cA(\cB)$ if $\cA$ is a basic class---of labelled objects obtained as follows. Choose an object $\gamma\in\cA$ to be the \emph{core} of the composition and let $k = |\gamma|$ be its size. Then pick a k-set of elements from $\cB$. Substitute each atom $v\in\gamma$ by an object $\gamma_v$ from the k-set, distributing distinct labels to the atoms of the composed object, i.e., the atoms in $\cup_{v\in\gamma}\gamma_v$. And forget the original labels on $\gamma$ and the $\gamma_v$. The composition construction is very powerful. For instance, it allows us to formulate the classical Set, Sequence, and Cycle constructions from basic classes. Indeed, the class of sequences (sets, cycles)  of objects in a class $\cA$ is simply the class $\Seq(\cA)$ ($\Set(\cA)$, $\Cyc(\cA)$, resp.). Sets, Cycles, and Sequences with
  a specific range for the number of components are also readily handled. We use the subscript notations $\Seq_{\geq k}(\cA)$, $\Set_{\geq k}(\cA)$, $\Cyc_{\geq k}(\cA)$, when the number of components is constrained to be at least some fixed value $k$.

\subsection{Counting series}\label{subsec:counting}
For labelled enumeration, the counting series is the exponential generating function, shortly the EGF, defined as
\begin{equation}
C(x):=\sum_n \frac{1}{n!}|\cC_n|x^n
\end{equation}
whereas for unlabelled enumeration the counting series is the \emph{ordinary generating function}, defined as
\begin{equation}
\widetilde{C}(x):=\sum_n |\widetilde{C}_n|x^n.
\end{equation}
In general, \emph{cycle index sums} are used for unlabelled enumeration as a convenient
refinement of ordinary generating functions. Cycle index sums are multivariate power series that preserve information on symmetries. A \emph{symmetry} of size $n$ on a class $\cC$ is a pair $(\sigma\in\mathfrak{S}_n,\gamma\in\cC_n)$ such that $\gamma$ is stable under the action of $\sigma$ 
(notice that $\sigma$ is allowed to be the identity). The corresponding weight is defined as $\prod_{i=1}^ns_i^{c_i}$, where $s_i$ is a formal variable and $c_i$ is the number of cycles of length $i$ in $\sigma$. The cycle index sum of $\cC$, denoted by $Z[\cC](s_1,s_2,\ldots)$, is the multivariate series defined as the sum of the weight-monomials over all symmetries on $\cC$. The ordinary generating function is obtained by substitution of $s_i$ by $x^i$: $$\widetilde{C}(x)=Z[\cC](x,x^2,\ldots).$$

\begin{figure}
\begin{tabular}[t]{|l|l|l|l|}
  \hline
  \textbf{Basic classes} & \textbf{Notation} & \textbf{EGF}  & \textbf{Cycle index sum} \\
  \hline
  \textbf{Neutral Class} & $\cC=1$ & $C(z) = 1$ & $Z[\cC] = 1$ \\
  \textbf{Atomic Class} & $\cC=\cZ$ & $C(z) = z$ & $Z[\cC] = s_1$ \\
  \textbf{Sequence} & $\cC=\Seq$ & $C(z) = \frac{1}{1 - z}$ & $Z[\cC] = \frac{1}{1-s_1}$ \\
  \textbf{Set} & $\cC=\Set$ & $C(z) = \exp (z)$ & $Z[\cC] =\exp\left( \sum_{r\geq 1} \frac{1}{r}s_r\right) $ \\
  \textbf{Cycle} & $\cC=\Cyc$ & $C(z) = \log\left(\frac{1}{1 - z}\right)$ & $Z[\cC] =\sum_{r\geq 1} \frac{\phi(r)}{r}\log\left(\frac{1}{1-s_r}\right) $ \\
  \hline
  \textbf{Construction} & \textbf{Notation} & \textbf{Rule for EGF} & \textbf{Rule for Cycle index sum}\\
  \hline
  \textbf{Union} & $\cC = \cA + \cB $ & $C(z) = A(z) + B(z)$ & $Z[\cC] = Z[\cA] + Z[\cB]$ \\
  \textbf{Product} & $\cC = \cA * \cB $ & $C(z) = A(z) \cdot B(z)$ & $Z[\cC] = Z[\cA] \times Z[\cB]$ \\
  \textbf{Composition} & $\cC = \cA \circ \cB $ & $C(z) = A(B(z))$ & $Z[\cC] = Z[\cA] \circ Z[\cB]$ \\
  \hline
\end{tabular}
\caption{Basic classes and constructions, with their translations to generating functions for labelled classes and to cycle index sums for unlabelled classes. For the composition construction, the notation $Z[\cA]\!\!\circ\!\! Z[\cB]$ refers to
the series {\small $Z[\cA]\!\circ\! Z[\cB](s_1,s_2,\ldots)=Z[\cA](Z[\cB](s_1,s_2,\ldots),Z[\cB](s_2,s_4,\ldots),Z[\cB](s_3,s_6,\ldots),\ldots)$}.}
\label{Basic}
\end{figure}

\subsection{Computation rules for the counting series}\label{sec:computation_rules}
The symbolic method provides for each basic class and each construction an explicit simple rule to compute the EGF (labelled enumeration) and the cycle index sum (unlabelled enumeration), as shown in Figure \ref{Basic}. These rules will allow us to convert our decomposition grammar into an enumerative strategy in an automatic way. As an example, consider  the  class $\cT$ of nonplane rooted trees.
Such a tree is made of a root vertex and a collection of subtrees pending from the root-vertex, which 
yields
$$
\cT=\cZ*\Set\circ\cT.
$$ 
For labelled enumeration, this is translated into the following equation satisfies by the EGF:
$$
T(x)=x\exp(T(x)).
$$
For unlabelled enumeration, this is translated into the following equation satisfied by the OGF
(via the computation rules for cycle index sums):
$$
\wt{T}(x)=x\exp\left(\sum_{r\geq 1}\frac1{r}\wt{T}(x^r)\right).
$$

In general, if a class $\cC$ is found to have a decomposition grammar, the rules of Figure~\ref{Basic} allow us to translate  the combinatorial description of the class into an equation-system satisfied by the counting series \emph{automatically} for both labelled and unlabelled structures. The purpose of this paper is to completely specify such a grammar to decompose any family of graphs into 3-connected components. 
Therefore we have to specify how the basic classes, constructions, and enumeration
tools have to be defined in the specific case of graph classes.

\subsection{Graph classes}
Let us first mention that the graphs we consider are allowed to have multiple edges
but no loops (multiple edges are allowed in the first formulation of the grammar, then we will explain how to adapt the grammar to simple graphs in 
Section~\ref{sec:simple}).
In the case of a class of graphs, we will need to take both vertices and edges into account. 
Accordingly, we consider a class  of graphs as a species of 
 combinatorial structures with two types of labelled atoms: vertices and edges. 
 In general we imagine  that if there are $n$ \emph{labelled} vertices and $m$ \emph{labelled} edges,
 then these labelled vertices carry distinct blue labels in $[1..n]$ and the edges carry distinct
 red labels in $[1..m]$~\footnote{If the graphs are simple, there is actually no need to label the edges, since two distinct edges are distinguished by the labels of their extremities.}. Hence, graph classes have to be treated in the extended
  framework of species with several types of atoms, see~\cite[Sec 2.4]{BeLaLe} 
  (we shortly review here how the basic 
  constructions and counting tools can be extended).
  
  For labelled enumeration  the exponential generating function (EGF) of a class of graphs is 
  $$
  G(x,y)=\sum_{n,m}\frac{1}{n!m!}|\cG_{n,m}|x^ny^m,
  $$
  where $\cG_{n,m}$ is the set of graphs in $\cG$ with $n$ vertices and $m$ edges.
For unlabelled enumeration (i.e., graphs are considered up  
  to relabelling the vertices  and the edges), the ordinary generating function (OGF) is
$$
  \wt{G}(x,y)=\sum_{n,m}|\wt{\cG}_{n,m}|x^ny^m,
  $$
  where $\wt{\cG}_{n,m}$ is the set of unlabelled graphs in the class that have 
  $n$ vertices and $m$ edges.
Cycle index sums can also be defined similarly as in the one-variable case, 
(as a sum of
weight-monomials) but the definition is more complicated, as well as the computation
rules, see~\cite{Wa82b}. In this article we restrict our attention to labelled enumeration
and postpone to future work the applications of our grammar to unlabelled enumeration.  
  
We distinguish three types of graphs: unrooted, vertex-pointed, and rooted.
 In an unrooted graph, all vertices and all edges are labelled.
 In a vertex-pointed graph, there is one distinguished vertex that is unlabelled, all
 the other vertices and edges are labelled.  In a rooted graph, 
 there is one distinguished edge---called the root---that is oriented, all
 the vertices are labelled except the extremities of the root,  and all edges are labelled except the root. A class of unrooted graphs is typically denoted by $\cG$, and the 
 associated vertex-pointed and rooted classes are respectively denoted $\cG'$ and $\ovr{\cG}$.
 Notice that $\cG'_{n,m}\simeq\cG_{n+1,m}$. The generating functions $G'$ of $\cG'$ and $\ovr{G}$ of $\ovr{\cG}$ satisfy:
 $$
 G'(x,y)=\partial_x G(x,y),\ \ \ \ \ \ovr{G}(x,y)=\frac{2}{x^2}\partial_yG(x,y).
 $$
A class of vertex-pointed graphs is called a \emph{vertex-pointed class}
and a class of rooted graphs is called a rooted class. In this article, all
vertex-pointed classes will be of the form $\cG'$, but we will consider rooted 
classes that are not of the form $\ovr{\cG}$; for such classes we require nevertheless
that the class is stable when reversing the direction of the root-edge.

 The basic graph classes are the following:
 \begin{itemize}
 \item
 The \emph{vertex-class} $v$ stands for the class made of 
 a unique graph that has a single vertex and no edge. The series is $(x,y)\mapsto x$. 
 \item
 The \emph{edge-class} $e$ stands for the class
  made of 
 a unique graph that has two unlabelled vertices connected by one directed labelled edge. The series is $(x,y)\mapsto y$. 
 \item
 The \emph{ring-class} $\cR$ stands for the
 class of \emph{ring-graphs}, which are cyclic chains of at least 3 edges. The
 series of $\cR$ is $(x,y)\mapsto \frac1{2}(-\log(1-xy)-xy-\frac1{2}\ \!x^2y^2)$. 
 \item
 The \emph{multi-edge-class}
 $\cM$ stands for the
 class of \emph{multi-edge graphs}, which consist of 2 labelled vertices connected
 by $k\geq 3$ edges. The
 series of $\cM$ is $(x,y)\!\mapsto\! \frac1{2}\ \!x^2(\exp(y)\!-\!1\!-\!y\!-\!\frac{y^2}{2})$.
 \end{itemize} 
 
 The constructions we consider for graph classes are the following: 
 \emph{disjoint union}, 
\emph{partitional product} (defined similarly as in the one-variable case), and now  two types of substitution:
 \begin{itemize}
 \item
 \emph{Vertex-substitution}: Given a graph class $\cA$ (which might be unrooted, vertex-pointed, or rooted) and a vertex-pointed class $\cB$, the class $\cC=\cA\circ_v\cB$
 is the class of graphs obtained by taking a graph $\gamma\in\cA$, called the core graph,
 and attaching at each labelled vertex $v\in\gamma$ a graph $\gamma_v\in\cB$,
 the vertex of attachment of $\gamma_v$ being the distinguished (unlabelled) vertex of $\gamma_v$.
We have 
 $$
 C(x,y)=A(xB(x,y),y),
 $$
 where $A$, $B$ and $C$ are respectively the exponential generating functions of $\cA$, $\cB$ and  $\cC$.
 \item
  \emph{Edge-substitution}: Given a graph class $\cA$ (which might be unrooted, vertex-pointed, or rooted) and a rooted class $\cB$, the class $\cC=\cA\circ_e\cB$
 is the class of graphs obtained by taking a graph $\gamma\in\cA$, called the core graph,
 and substituting each labelled edge $e=\{u,v\}$ (which is implicitly given an orientation)
  of $\gamma$ by a graph $\gamma_e\in\cB$, thereby identifying the origin of the root of  
  $\gamma_e$ with $u$ and the end of the root of $\gamma_e$ with $v$.
  After the identification, the root edge of $\gamma_e$ is deleted.
We have 
 $$
 C(x,y)=A(x,B(x,y)),
 $$
 where $A$, $B$ and $C$ are respectively the generating functions of $\cA$, $\cB$ and $\cC$.
 \end{itemize}

\section{Tree decomposition and dissymmetry theorem}\label{DissymThms}

The dissymmetry theorem for trees \cite{BeLaLe} makes it possible to express the class of unrooted trees in terms of classes of rooted trees. Precisely, 
let $\cA$ be the class of tree, and let us define the following associated rooted families: $\cA_{\circ}$ is the class of trees
where a node is marked,  $\cA_{\circ-\circ}$ is the class of trees where an edge is marked, and $\cA_{\circ\rightarrow\circ}$ is the class
of trees where an edge is marked and is given a direction.   
Then the class $\cA$ is related to these three associated rooted classes by the following identity:
\begin{equation}\label{DissymTrees}
\cA + \cA_{\circ\rightarrow\circ} \simeq \cA_{\circ} + \cA_{\circ-\circ}.
\end{equation}

 The theorem is named after the dissymmetry resulting in a tree rooted anywhere other than at its centre, see~\cite{BeLaLe}.
Equation (\ref{DissymTrees}) is an elegant and flexible counterpart to the dissimilarity equation discovered by Otter~\cite{Ot48}; as we state in Theorem \ref{DissymTreesDecomp} below, it can easily be extended to classes for which a tree can be \emph{associated} with each object in the class.

A \emph{tree-decomposable} class is a class $\cC$ such that to each object $\gamma \in \cC$ is associated a tree $\tau(\gamma)$ whose nodes are distinguishable in some way (e.g., using the labels on the vertices of $\gamma$). Denote by $\cC_{\circ}$ the class of objects of $\cC$ where a node of $\tau(\gamma)$ is distinguished, by $\cC_{\circ-\circ}$ the class of objects of $\cC$ where an edge of $\tau(\gamma)$ is distinguished, and by $\cC_{\circ\rightarrow\circ}$ the class of objects of $\cC$ where an edge of $\tau(\gamma)$ is distinguished and given
a direction. The principles and proof of the dissymmetry theorem can be straightforwardly extended to any tree-decomposable class, giving rise to the 
following statement.

\begin{theorem}\label{DissymTreesDecomp}(Dissymmetry theorem for tree-decomposable classes) Let $\cC$ be a tree-decomposable class. Then
\begin{equation}\label{eq:diss_general}
\cC + \cC_{\circ\rightarrow\circ} \simeq \cC_{\circ} + \cC_{\circ-\circ}.
\end{equation}
\end{theorem}

Note that, if the trees associated to the graphs in $\cC$ are bipartite, then  
$\cC_{\circ\rightarrow\circ}\simeq 2\cC_{\circ-\circ}$. Hence, Equation~\eqref{eq:diss_general} simplifies to
\begin{equation}
\cC \simeq \cC_{\circ} - \cC_{\circ-\circ}.
\end{equation}
(At the upper level of generating functions, this reflects the property that the number of vertices in a tree exceeds the number of edges by one.) 

\section{Tutte's decomposition and beyond: Decomposing a graph into 3-connected components}\label{DeGr3-cComps}
In this section we recall Tutte's decomposition~\cite{Tutte} 
of a graph into 3-connected components,
which we will translate into a grammar in Section~\ref{sec:Decomp Grammar}. 
The decomposition works in three levels: (i) standard decomposition of a graph into connected 
components, (ii) decomposition of a connected graph into 2-connected blocks that are
articulated around vertices, (iii)
decomposition of a 2-connected into 3-connected components that are articulated around 
(virtual) edges. 

A nice feature of Tutte's decomposition is that the second and third level are ``tree-like''
decompositions, meaning that the ``backbone'' of the decomposition is a tree.
The tree associated with (ii) is called the Bv-tree, and the tree associated with (iii)
is called the RMT-tree (the trees are named after the possible types of the nodes).
The tree-property of the decompositions will
 enable us to apply the dissymmetry theorem---Theorem~\ref{DissymTreesDecomp}---in order to write down the grammar.
As we will see in Section~\ref{grammarCfrom2C}, 
writing the grammar will require the canonical decomposition of \emph{vertex-pointed}
2-connected graphs. It turns out that a smaller backbone-tree (smaller than for \emph{unrooted} 2-connected graphs) is more convenient in order to
apply the dissymmetry theorem, thereby simplifying  the decomposition process for vertex-pointed 2-connected graphs.
Thus in Section~\ref{sec:root_RMT_tree} we introduce these smaller trees, called restricted RMT-trees (to our knowledge, these trees have not been considered before).

\subsection{Graphs and connectivity}\label{subsec:concepts}
We give here a few definitions on graphs and connectivity, 
following Tutte's terminology~\cite{Tutte}.
The vertex-set (edge-set) of a graph $G$ is denoted by $V(G)$ ($E(G)$, resp.). A subgraph of a graph $G$ is a graph $G'$ such that $V(G')\subset V(G)$, $E(G')\subset E(G)$, and any vertex incident to an edge in $E(G')$ is in $V(G')$. Given an edge-subset $E'\subset E(G)$, the corresponding \emph{induced graph} is the subgraph $G'$ of $G$ such that $E(G')=E'$ and $V(G')$ is the set of vertices incident to edges in $E'$; the induced graph is denoted by $G[E']$.

A graph is \emph{connected} if any two of its vertices are connected by a path. A \emph{1-separator} of a graph $G$ is given by a partition of $E(G)$ into two nonempty sets $E_1,E_2$ such that $G[E_1]$ and $G[E_2]$ intersect at a unique vertex $v$; such a vertex is called \emph{separating}. A graph is \emph{2-connected} if it has at least
two vertices and no 1-separator. 
Equivalently (since we do not allow any loop), a  2-connected graph $G$
 has at least two vertices and the deletion of any vertex does not disconnect $G$.
A \emph{2-separator} of a graph is given by a partition of $E(G)$ into two subsets $E_1,E_2$ each of cardinality at least 2, such that $G[E_1]$ and $G[E_2]$ intersect at two vertices $u$ and $v$; such a pair $\{u,v\}$ is called a \emph{separating vertex pair}. A graph is \emph{3-connected} if it has no 2-separator and has at least 4 vertices. (The latter condition is convenient for our purpose, as it 
 prevents any ring-graph or multiedge-graph from being 3-connected.) Equivalently, a  3-connected graph $G$ has at least 4 vertices, no loop nor multiple edges, and the deletion of any two vertices does not disconnect $G$.

\subsection{Decomposing a connected graph into 2-connected ones}\label{subsec:1to2}
There is a well-known decomposition of a graph into 2-connected components, which 
is described in several books~\cite{Ha,Diestel,Mohar,Tutte}.

\begin{figure}\label{Bv-tree}
\begin{center}
\includegraphics[width=14cm]{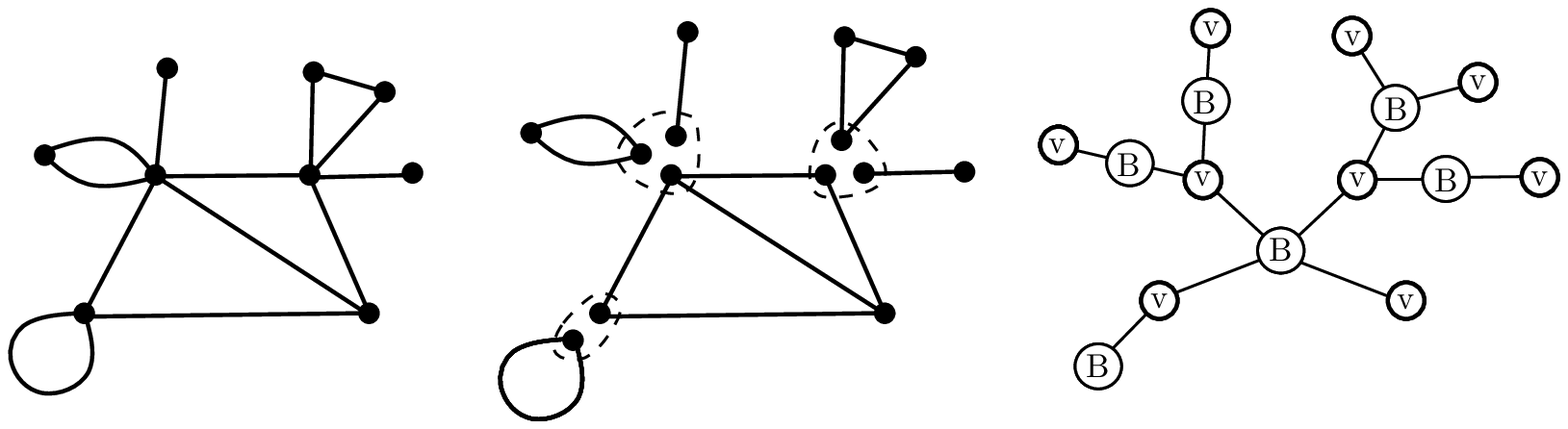}
 \end{center}
 \caption{Decomposition of a connected graph into blocks, and the associated Bv-tree.}
 \label{fig:blocks}
 \end{figure}

Given a connected graph $C$, a \emph{block} of $C$ is a maximal 2-connected induced subgraph of $C$. The set of blocks of $C$ is denoted by $\mathfrak{B}(C)$. 
A vertex $v\in C$ is said to be \emph{incident} to a block $B\in\mathfrak{B}(C)$ if $v$ belongs to $B$. The \emph{Bv-tree} of $C$ describes
 the incidences between vertices and blocks of $C$, i.e., 
it is a bipartite graph $\tau(C)$ with  node-set $V(C)\cup\mathfrak{B}(C)$, 
and edge-set given by the incidences between the vertices and the blocks of $C$, see Figure~\ref{fig:blocks}. The graph $\tau(C)$ is actually 
a tree, as shown for instance in \cite{Tutte,Mohar}. 
Conversely, take a collection $\mathfrak{B}$ of 2-connected graphs, called blocks, and a vertex-set $V$ such that every vertex in $V$ is in at least one block and the graph of incidences between blocks and vertices is a tree $\tau$. Then the resulting  graph is connected and has $\tau$ as its Bv-tree.  Consequently, connected graphs can be \emph{identified} with their tree-decompositions into blocks,  
which  will be very useful for deriving decomposition grammars.

\subsection{Decomposing a 2-connected graph into 3-connected ones}\label{subsec:2to3}
In this section we recall Tutte's decomposition of a 2-connected graph into
3-connected components~\cite{Tu63}.
 A similar decomposition has also been described by Hopcroft and 
 Tarjan~\cite{HoTa73}, 
 however they use a split-and-remerge process, whereas  Tutte's method
 only involves (more restrictive) split operations. We follow here the 
 presentation of Tutte.

First, one has to define connectivity modulo a pair of vertices. Let $G$ be a 2-connected graph  and $\{u,v\}$ a pair of vertices of $G$. Then $G$ is said to be \emph{connected modulo $[u,v]$} if there exists no partition of $E(G)$ into two nonempty sets $E_1,E_2$ such that $G[E_1]$ and $G[E_2]$ intersect only at $u$ and $v$. Being non-connected modulo $[u,v]$ means either that $u$ and $v$ are adjacent or that the deletion of  $u$ and $v$ disconnects the graph.

\begin{figure}\label{split}
\begin{center}
\includegraphics[width=10cm]{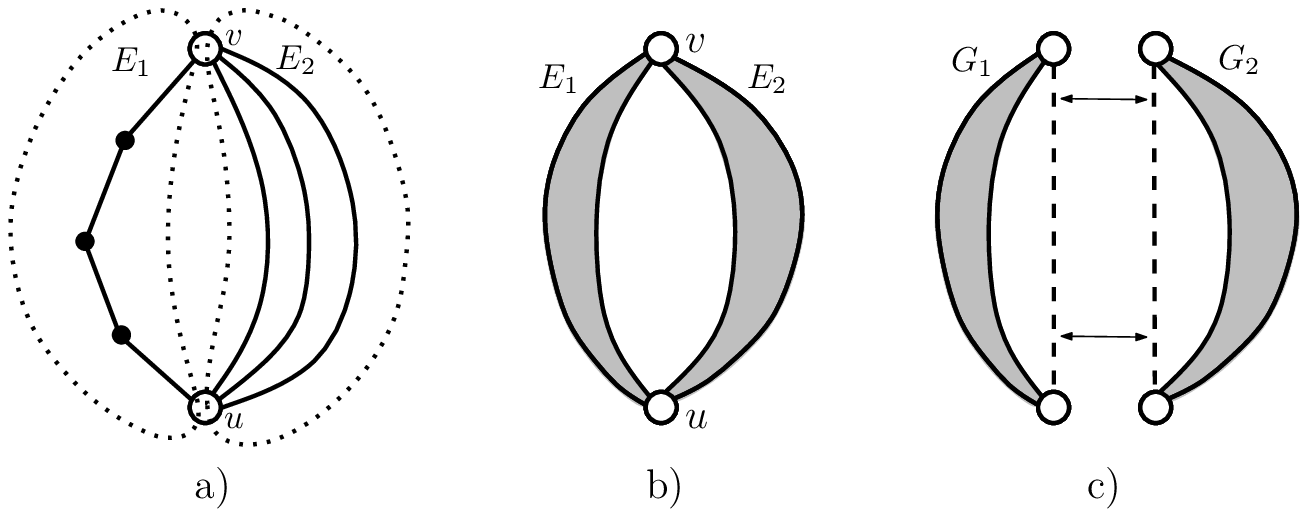}
 \end{center}
 \caption{(a) Example of a split candidate. (b) Splitting the edge-set. (c) Splitting a graph along a virtual edge.}
 \label{fig:split}
 \end{figure}

Consider a 2-separator $E_1,E_2$ of a 2-connected graph $G$, with $u,v$ the corresponding separating vertex-pair. Then $E_1,E_2$ is called a \emph{split-candidate}, denoted by $\{E_1,E_2,u,v\}$,  if $G[E_1]$ is connected modulo $[u,v]$ and $G[E_2]$ is 2-connected.  Figure~\ref{fig:split}(a) gives an example of a split-candidate, where $G[E_1]$ is connected modulo $[u,v]$ but not 2-connected, while $G[E_2]$ is 2-connected but not connected modulo $[u,v]$.

As described below, split candidates make it possible to completely decompose a 2-connected graph into 3-connected components. We consider here only 2-connected graphs with at least 3 edges (graphs with less edges are degenerated for this decomposition). Given a split candidate $S=\{E_1,E_2,u,v\}$ in a 2-connected graph $G$  (see Figure~\ref{fig:split}(b)),
the corresponding \emph{split operation} is defined as follows, see Figure~\ref{fig:split}(b)-(c): 
\begin{itemize}
\item
an edge $e$, called a \emph{virtual edge}, is added between $u$ and $v$, 
\item
the graph $G[E_1]$ is separated from the graph $G[E_2]$ 
by cutting along the edge $e$. 
\end{itemize}
Such a split operation yields two graphs $G_1$ and $G_2$, see Figure~\ref{fig:split}(d), which correspond respectively to $G[E_1]$ and $G[E_2]$ together with $e$ as a real edge. The graphs $G_1$ and $G_2$ are said to be \emph{matched} by the virtual edge $e$. It is easily checked that $G_1$ and $G_2$ are 2-connected (and have at least 3 edges). The splitting process can be repeated until no split candidate remains left.

\begin{figure}
\begin{center}
\includegraphics[width=16cm]{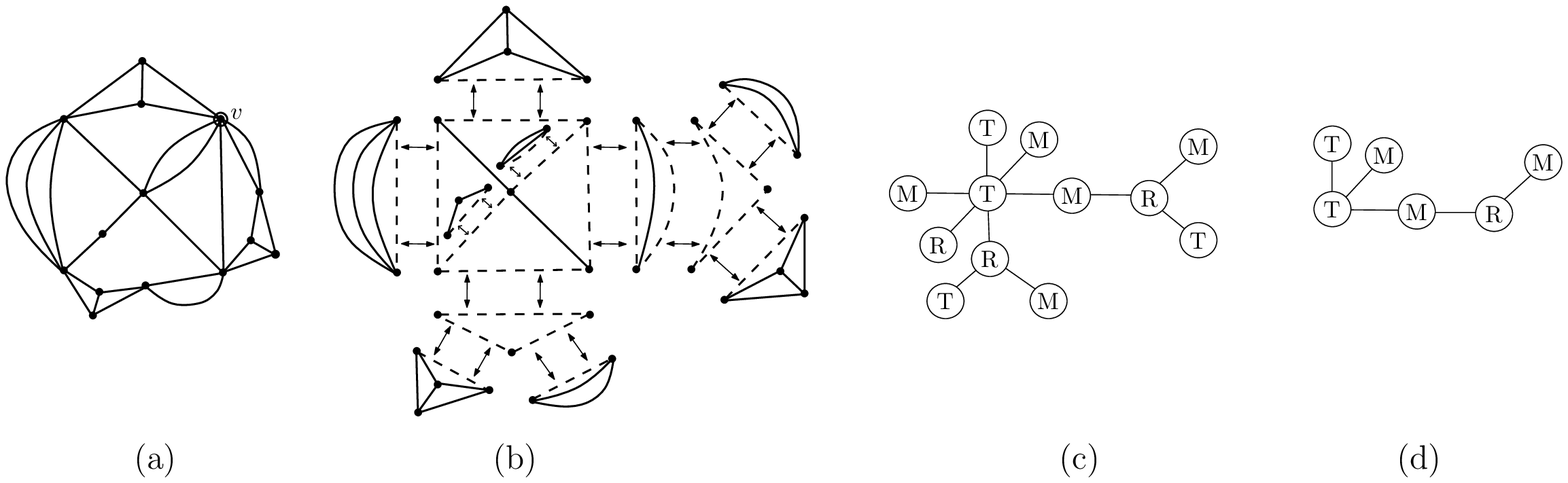}
 \end{center}
 \caption{(a) A 2-connected graph, (b) decomposed into bricks. (c) The associated RMT-tree. (d) The associated restricted RMT-tree if the graph
 is pointed at $v$.}
\label{fig:bricks}
 \end{figure}

As shown by Tutte in~\cite{Tutte}, the structure resulting from the split operations is independent of the order in which they are performed. It is a collection of graphs, called the \emph{bricks} of $G$, which 
are articulated around virtual edges, see Figure~\ref{fig:bricks}(b). 
By definition of the decomposition, each brick has no split candidate; 
Tutte has shown that such graphs are either multiedge-graphs (M-bricks) or ring-graphs (R-bricks), or 3-connected graphs with at least 4 vertices (T-bricks).

The \emph{RMT-tree} of $G$ is the graph $\tau(G)$ whose nodes are the bricks of $G$ and whose edges correspond to the virtual edges of $G$ (each virtual edge matches two bricks),
see Figure~\ref{fig:bricks}. 
The graph $\tau(G)$ is indeed a tree~\cite{Tutte}. 
By maximality of the decomposition, it is easily checked that $\tau(G)$ has no two R-bricks adjacent nor two M-bricks adjacent. 

Call a \emph{brick-graph} a graph that is either a ring-graph or a multi-edge graph or a 3-connected graph, with again the letter-triple $\{R,M,T\}$ to refer to the type of the brick. The inverse process of the split decomposition consists in taking a collection of brick-graphs and a collection of edges, called virtual edges, so that each virtual edge belongs to two bricks, and so that the graph $\tau$ with vertex-set the bricks and edge-set  the virtual edges (each virtual edge matches two bricks) is a tree avoiding two R-bricks or  two M-bricks being adjacent. Then the resulting graph, obtained by matching the bricks along  virtual edges and then erasing the virtual edges, is a 2-connected graph that has $\tau$ as its RMT-tree. Hence, 2-connected graphs with at least 3 edges can be 
identified with their RMT-tree, which again will be useful for writing down a decomposition grammar.

\subsection{The restricted RMT-tree}\label{sec:root_RMT_tree}
The grammar to be written in Section \ref{sec:Decomp Grammar} requires to decompose not only
unrooted 2-connected graphs, but also \emph{vertex-pointed} 2-connected graphs. 
It turns out that these vertex-pointed 2-connected graphs are much more convenient to decompose using a subtree of the RMT-tree.

The \emph{restricted RMT-tree} of a 2-connected graph $G$ (with at least 3 edges) rooted at a vertex $v$, is defined as the subgraph $\tau'(G)$ of the RMT-tree $\tau(G)$ induced by the bricks containing $v$ and by the edges of $\tau(G)$ connecting two such bricks.

\begin{lemma}\label{rootRMT}
The restricted RMT-tree of a vertex-pointed 2-connected graph with at least 3 edges is a tree.
\end{lemma}
\begin{proof}
Let $G$ be a vertex-pointed 2-connected graph with at least 3 edges. Let $\tau(G)$ be the RMT-tree of $G$ and let $\tau'(G)$ be the restricted RMT-tree of $G$. The pointed vertex is denoted by $v$. As $\tau'(G)$ is a subgraph of the tree $\tau(G)$, it is enough to show that $\tau'(G)$ is connected for it to be a tree. Recall that a virtual edge $e$ corresponds to splitting $G$ into two graphs $G_1=G[E_1]+e$ and $G_2=G[E_2]+e$, where $E_1,E_2$ is a 2-separator of $G$. The two subtrees $T_1$ and $T_2$ attached at each extremity of the virtual edge correspond to the split-decomposition of $G_1$ and $G_2$, respectively. Hence, the pointed vertex 
$v$, if not  incident to the virtual edge $e$, is either a vertex of $G_1\backslash e$ or is a vertex of $G_2\backslash e$. In the first (second) case, $\tau'(G)$ is contained in $T_1$ ($T_2$, respectively). Hence, if an edge of $\tau(G)$ is not in $\tau'(G)$, then $\tau'(G)$ does not overlap simultaneously with the two subtrees attached at each extremity of that edge. This  property ensures that $\tau'(G)$ is connected.
\end{proof}

Having proved that the restricted RMT-tree is indeed a tree and not a forest, we will be able to use the dissymmetry theorem---Theorem \ref{DissymTreesDecomp}---in order to write a decomposition grammar for the class of vertex-pointed 2-connected graphs. 
The restricted RMT-tree turns out to be much better adapted for this purpose than the RMT-tree.

\section{Decomposition Grammar}\label{sec:Decomp Grammar}
In this section we translate Tutte's decomposition into an explicit grammar. Thanks to this grammar,
counting a family of graphs reduces to counting the 3-connected subfamily, which turns
out to be a fruitful strategy in many cases, in particular for planar graphs, as
we will see Section~\ref{sec:planar}.

Given a graph family $\cG$, our grammar corresponds at the first level to the connected components,
at the second level to the decomposition of a connected graph into 2-connected blocks, and at the third level to the decomposition of a 2-connected graph into 3-connected components. The first level is classic,
the second level already makes use of the dissymmetry theorem, it 
 is implicitly used by Robinson \cite{Ro70}, 
 and appears explicitly in the work by Leroux \cite{Le88, BeLaLe}. 
 The third level is new (though Leroux et al~\cite{GLLW08} have recently independently derived general equation
 systems relating the series of 2-connected graphs and 3-connected graphs of a given class).
 As we will see, it makes an even more extensive use of the dissymmetry theorem than the second level.

 We define the following subfamilies of $\cG$:
 \begin{itemize}
 \item
 The class $\cG_1$ is the subfamily of graphs in $\cG$ that are
 connected and have at least one vertex.
 \item
 The class $\cG_2$ is the subfamily of graphs in $\cG$ that are
 2-connected and have at least two vertices. Multiple edges are allowed.
(The smallest possible such graph is the link-graph that has two vertices connected by one edge.)
\item
 The class $\cG_3$ is the subfamily of graphs in $\cG$ that are
 3-connected and have at least four vertices. 
 (The smallest possible such graph is the tetrahedron.)
 \end{itemize}

A class $\cG$ of graphs is said to be \emph{stable under Tutte's decomposition} if it satisfies the following 
 property: 
 \begin{center}
 ``any graph $G$ is in $\cG$ iff all 3-connected components of $G$ are in $\cG$''.
  \end{center}
 Notice that a class of graphs stable under Tutte's decomposition satisfies the following
 properties:
 \begin{itemize}
 \item
 a graph $G$ is in $\cG$ iff all its connected components are in $\cG_1$,
 \item
 a graph $G$ is in $\cG_1$ iff all its 2-connected components are in $\cG_2$,
 \item
a graph $G$ is in $\cG_2$ iff all its 3-connected components are in $\cG_3$.
 \end{itemize}

\subsection{General from connected graphs}\label{subsec:generalfromconnected}
The first level of the grammar is classic.
A graph is simply the collection of its connected components, which translates to:

\begin{equation}\label{GfromC}
\cG = \Set(\cG_1).
\end{equation}

\subsection{Connected from 2-connected graphs}\label{grammarCfrom2C}
In order to write down the second level, i.e., decompose the connected class $\cC:=\cG_1$, 
we define the following classes:  $\cC_{B}$ is the class of graphs in $\cG_1$ with a distinguished block, $\cC_{v}$  is the class of graphs in $\cG_1$ with a distinguished vertex, and $\cC_{{Bv}}$  is the class of graphs in $\cG_1$ with a distinguished incidence block-vertex. In other words, $\cC_{B}$, $\cC_v$, and $\cC_{Bv}$ correspond to graphs in $\cG_1$ where one distinguishes  in the associated Bv-tree, respectively, a v-node, a B-node,
and an edge. The generalized dissymmetry theorem yields the following relation between $\cC=\cG_1$ and the auxiliary rooted classes:
$$
\cC + \cC_{Bv} = \cC_v + \cC_B,
$$
which can be rewritten as 
\begin{equation}\label{eq:G1}
\cC = \cC_v + \cC_B-  \cC_{Bv}.
\end{equation}

Clearly the class $\cC'$ is related to $\cC_{v}$ by $\cC_v = v * \cC'$.
To decompose $\cC'$, we observe that the pointed vertex gives a starting point for
 a recursive decomposition.
Precisely, from the block decomposition described in Section~\ref{subsec:1to2}, 
any vertex-pointed connected graph is obtained as follows:
take a collection of vertex-pointed 2-connected graphs attached together at their marked vertices,
 and attach a vertex-pointed connected graph at each non-pointed vertex of these 2-connected graphs.
 (Clearly the 2-connected graphs correspond to the blocks incident to the pointed vertex in the resulting
 graph.)

 This recursive decomposition translates to the equation
\begin{equation}\label{eq:G1p}
\cC' = \Set ( \cGbp\circ_{v} \cC').
\end{equation}

Similarly, each graph in $\cC_{B}$ is obtained in a unique way by taking a block in $\cG_2$ and attaching at each vertex of the block a vertex-pointed connected graph in $\cC'$, which yields
 \begin{equation}\label{eq:G1B}
\cC_{B} = \cG_2 \circ_{v} \cC'.
 \end{equation}

Finally, each graph in $\cC_{Bv}$ is obtained from a vertex-pointed block in $\cGbp$ by attaching at  each vertex of the block ---even the root vertex--- a vertex-pointed connected graph, which yields
\begin{equation}\label{eq:G1Bp}
\cC_{Bv} = (v*\cGbp)\circ_{v}\cC'.
\end{equation}
The grammar to decompose a class of connected graphs into 2-connected components results
from the concatenation of Equations~\eqref{eq:G1},~\eqref{eq:G1p},~\eqref{eq:G1B}, and~\eqref{eq:G1Bp}.

A similar grammar is given in the book of Bergeron, Labelle and Leroux \cite{BeLaLe}. Notice that there
are two terminal classes in this grammar, the class $\cG_2$ and the class $\cGbp$.

\subsection{2-Connected from 3-connected graphs}\label{grammar2Cfrom3C}
In this section we start to describe the new contributions of this article,
namely the decomposition grammars for $\cG_2$ and $\cGbp$. 

Let us begin with $\cG_2$.
Again we have to define auxiliary classes that correspond to the different ways to distinguish a node
or an edge in the RMT-tree.
Let $\cBthree$ be the class of graphs in $\cG_2$ with at least 3 edges
(i.e., those whose RMT-tree is not empty). Since we consider graph classes stable under Tutte's decomposition, the link-graph $\ell_1$ and
the double-link graph $\ell_2$ (which have counting series $x^2y/2$ and $x^2y^2/2$, respectively) 
are in $\cG_2$, hence
\begin{equation}\label{eq:Bthree}
\cG_2=\ell_1+\ell_2+\cBthree.
\end{equation}

Next we decompose $\cBthree$ using the RMT-tree.
Let $\cBv$ ($\cBe$, $\cBh$) be the class of graphs in $\cBthree$ such that the RMT-tree carries a distinguished node 
(edge, directed edge, resp.). Theorem \ref{DissymTreesDecomp} yields
\begin{equation}\label{eq:diss_2conn}
\cBthree =  \cBv+\cBe-\cBh.
\end{equation}
The class $\cBv$ is naturally partitioned into 3 classes $\cBR$, $\cBM$, and $\cBT$, depending on the type of the distinguished node (R-node,  M-node, or T-node). 
Similarly, the class $\cBe$ is partitioned into 4 classes $\cBRM$, $\cBRT$, $\cBMT$, and $\cBTT$ 
(recall that a RMT-tree has no two adjacent R-bricks nor two adjacent M-bricks); 
and $\cBh$ is partitioned into 7 classes $\cBRtoM$, $\cBMtoR$, $\cBRtoT$, $\cBTtoR$, $\cBMtoT$,
$\cBTtoM$, and $\cBTtoT$. Notice that $\cBRtoM\simeq\cBMtoR\simeq\cBRM$, $\cBRtoT\simeq\cBTtoR\simeq\cBRT$, and $\cBMtoT\simeq\cBTtoM\simeq\cBMT$. Hence, Equation~(\ref{eq:diss_2conn}) 
is rewritten as
\begin{equation}\label{DissymmThreeGraphsEq}
\cBthree =  \cBR + \cBM + \cBT -\cBRM - \cBRT - \cBMT-\cBTtoT +\cBTT.
\end{equation}

\subsubsection{Networks}
In order to decompose the classes on the right-hand-side 
of Equation~(\ref{DissymmThreeGraphsEq}), we 
first have to decompose the class  of rooted 2-connected graphs in $\cG$, more
precisely we need to specify a grammar for a class of objects closely related to $\ovr{\cG_2}$,
which are called networks.
A network is defined as a connected graph arising from a graph in $\ovr{\cG_2}$  
by deleting the root-edge; the origin and end of the root-edge are respectively called the $0$-pole and the 
  $\infty$-pole of the network. 
The associated class  is classically denoted by $\cD$
in the literature~\cite{Wa82a}.
Observe that the only rooted 2-connected graph disconnected by root-edge deletion 
 is the rooted link-graph.
Hence
$$
\cGbr=1+\cD,
$$
where the rooted link-graph has weight $1$ instead of $e$ because, in a rooted class, the rooted edge
is considered as unlabelled, i.e.,  is not counted in the size parameters.
(We will see in Section~\ref{sec:simple} that the link between $\cD$ and $\cGbr$ is a bit more complicated if 
multiple edges are forbidden.)

As discovered by Tracktenbrot~\cite{trak} a few year's before Tutte's book appeared, 
the class of networks with at least 2 edges (recall that the root-edge has been deleted) is naturally partitioned into 3 subclasses: $\cS$ for \emph{series networks}, $\cP$ for \emph{parallel networks}, and $\cH$ for \emph{polyhedral networks}:
\begin{equation}\label{eq:D}
\cD=e+\cS+\cP+\cH.
\end{equation}
With our terminology of RMT-tree, the three situations correspond to 
  the root-edge belonging to a R-brick, M-brick, or 
  T-brick, respectively~\footnote{Actually Trackhtenbrot's 
  decomposition can be seen as Tutte's decomposition restricted to rooted 2-connected graphs.}. 
    In a similar way as for the class $\cGcp$ in Section~\ref{grammarCfrom2C}, 
    the root-edge gives a starting
  point for a recursive decomposition. Clearly, as there is no edge R-R in the RMT-tree,  
  each series network is obtained as a collection of at least two non-series networks connected as a chain (the $\infty$-pole of a network is identified with the $0$-pole of the following network in the chain):
    \begin{equation}\label{eq:S}
  \cS=(\cD-\cS)* v*\cD.
  \end{equation}
  Similarly, as there is no edge M-M in the RMT-tree, 
  each parallel-network is obtained as a collection of at least two non-parallel networks sharing the same $0$- and $\infty$-poles:
    \begin{equation}\label{eq:P}
    \cP=\Set_{\geq 2}(\cD-\cP).
    \end{equation}
  Finally, each polyhedral network is obtained as a rooted 3-connected graph 
  where each non-root edge 
  is substituted by a network, which yields:
  \begin{equation}\label{eq:H}
  \cH=\ovr{\cG_3}\circ_e\cD.
  \end{equation}
  The resulting decomposition grammar for $\cD$ is obtained as the concatenation of Equations~\eqref{eq:D},~\eqref{eq:S},~\eqref{eq:P}, and~\eqref{eq:H}.  
 This grammar has been known since Walsh \cite{Wa82a}.
 Notice that the only terminal class is the 3-connected class $\ovr{\cG_3}$.


\subsubsection{Unrooted 2-connected graphs}
We can now  specify the decompositions of the families on the right-hand-side 
of~\eqref{DissymmThreeGraphsEq}.
 Recall that $\cR$ is the class of ring-graphs (polygons) and  $\cM$ is the class of 
  multiedge graphs with at least 3 edges. Given a graph in $\cBR$, each edge $e$ of the distinguished R-brick is either a real edge or a virtual edge; 
  in the latter case the graph attached on the other side of $e$ (i.e., the side not
incident to the rooted R-brick) is naturally rooted at $e$; it is thus identified with a network (upon choosing an orientation of $e$),
precisely it is a non-series network, as there are no two R-bricks adjacent. Hence 
\begin{equation}
\cBR=\cR\circ_e(\cD-\cS).
\end{equation}
Similarly we obtain
\begin{equation}
\cBM=\cM\circ_e(\cD-\cP)=(v^2*\Set_{\geq 3}(\cD-\cP))/\bullet\leftrightarrows\bullet,
\end{equation}
where the last notation means ``up to exchanging the two pole-vertices of the multiedge component'', and  
\begin{equation}
\cBT=\cG_3\circ_e\cD,
\end{equation}
to be compared with $\cH=\ovr{\cG_3}\circ_e\cD$. Next we decompose 2-connected graphs with a distinguished edge in the  RMT-tree. 
Consider the class $\cBRM$. The distinguished edge of the RMT-tree corresponds to a virtual edge $\{u,v\}$ matching one R-brick and one M-brick, such that the two bricks are attached at $\{u,v\}$. Upon fixing an orientation of the virtual edge $\{u,v\}$, there is a series-network on one side of $\{u,v\}$ and a parallel-network on the other side. Notice that such a construction has to be considered up to 
 orienting $\{u,v\}$, i.e., up to 
exchanging 
the two poles $u$ and $v$ (notation $/\bullet\leftrightarrows\bullet$). We obtain
\begin{equation}
\qquad\cBRM = (\mathcal{S} * \mathcal{P})/\bullet\leftrightarrows\bullet.
\end{equation}

\vspace{-.2cm}

\noindent Similarly,

\vspace{-.7cm}

 $$
 \cBMT=(\cP*\cH)/\bullet\leftrightarrows\bullet,\ \ \cBRT=(\cS*\cH)/\bullet\leftrightarrows\bullet,\ \ \cBTtoT=(\cH*\cH)/\bullet\leftrightarrows\bullet,
 $$
 
 \vspace{-.6cm}
 
 $$
\ \ \cBTT = (\mathcal{H} * \mathcal{H})/(\bullet\leftrightarrows\bullet , H\leftrightarrows H),
 $$
 
 \vspace{-.1cm}
 
\noindent where the very last notation means ``up to orienting the distinguished virtual edge''.

\subsubsection{Vertex-pointed 2-connected graphs}\label{sec:grammarPointed2conn}
Now we decompose the class $\cGbp$ of vertex-pointed 2-connected graphs 
(recall that $\cG_2\,\!\!'$ is, with $\cG_2$, one of the two terminal classes in the decomposition grammar for connected graphs into 2-connected components). We proceed in a similar way as for $\cG_2$, with the important difference that we use the restricted RMT-tree instead of the RMT-tree. Observe that the class of graphs in $\cG_2\ \!\!\!'$ with at 
least 3 edges (those whose RMT-tree is not empty) is the derived class $\cV:=\cB'$. By deriving the identity~\eqref{eq:Bthree}, we get
\begin{equation}
\cGbp=\ell_1\ \!\!\!'+\ell_2\ \!\!\!'+\cV=v*e+v*\Set_2(e)+\cV.
\end{equation} 
We denote by $\cV_{\circ}$ ($\cV_{\circ-\circ}$, $\cV_{\circ\rightarrow\circ}$) the class  of graphs in $\cGbp$ where the associated restricted RMT-tree
carries a distinguished 
node (edge, oriented edge, resp.). Theorem~\ref{DissymTreesDecomp} yields
\begin{equation}
\cV=\cV_{\circ}+\cV_{\circ-\circ}-\cV_{\circ\rightarrow\circ}.
\end{equation}

Notice that $\cV_{\circ}\neq(\cB_{\circ})'$, $\cV_{\circ-\circ}\neq(\cB_{\circ-\circ})'$, and $\cV_{\circ\rightarrow\circ}\neq(\cB_{\circ\rightarrow\circ})'$, due to the fact that the associated tree is not the same for $\cV=\cB'$ and for $\cB$ (actually, taking the derivative of the grammar for $\cB$ would produce a 
much more complicated grammar for $\cB'$ than the one we will obtain using the restricted RMT-tree). We partition the classes 
$\cV_{\circ}$, $\cV_{\circ-\circ}$, and $\cV_{\circ\rightarrow\circ}$ according to the 
types of the bricks incident to the root, and proceed with a decomposition in each case; the arguments are
very similar as for the decomposition of $\cG_2$. Take the example of $\cV_T$. Since the pointed vertex of the graph is incident to the marked T-brick (this is 
where it is very nice to consider the restricted RMT-tree instead of the RMT-tree), we have $$\cV_T=\cG_3\ \!\!\!'\circ_e\cD,$$ 
to be compared with $\cB_T=\cG_3\circ_e\cD$. Similarly we obtain

\vspace{-.4cm}

$$
\cV_R=\cR'\circ_e(\cD-\cS),\ \ \ \ \ \ \ \  \cV_M=\cM'\circ_e(\cD-\cP)=v*\Set_{\geq 3}(\cD-\cP),
$$

\vspace{-.6cm}

$$\cV_{R-M} = v*\mathcal{S} * \mathcal{P},\qquad\cV_{R-T} = v*\mathcal{S} * \mathcal{H},\qquad\cV_{M-T}=v*\mathcal{P} * \mathcal{H},$$

\vspace{-.6cm}

$$\qquad\cV_{T\rightarrow T}=v*\mathcal{H} *\mathcal{H},\qquad\cV_{T-T}=(v*\mathcal{H} * \mathcal{H})/H\leftrightarrows H.$$

\begin{figure}
\begin{center}
\includegraphics[width=8cm]{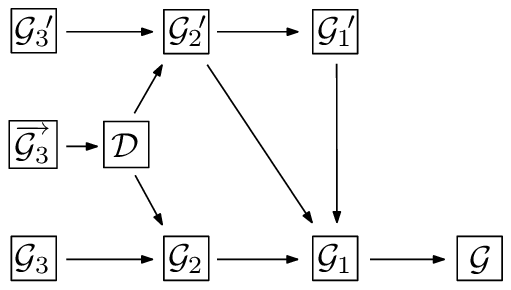}
\end{center}
\caption{The main dependencies in the grammar.}
\label{fig:deps}
\end{figure}

We have now a generic grammar to decompose a family of graphs into  3-connected components;  the complete grammar is shown in Figure~\ref{GrammarFig}, and the main dependencies are shown in Figure~\ref{fig:deps}. 
 Observe that, except for the 
easy basic classes $v$, $e$, $\ell_1$, $\ell_2$, $\Set_{\geq k}$, $\cR$, $\cM$, $\cR'$, $\cM'$, the only terminal classes are the 3-connected classes $\cG_3$, $\cG_3'$, 
and $\ovr{\cG_3}$.

\subsection{Adapting the grammar for families of simple graphs}\label{sec:simple}
The grammar has been described for a family $\cG$ satisfying the stability condition under Tutte's decomposition, and where multi-edges are allowed.
It is actually very easy to adapt the grammar for the corresponding \emph{simple} family of graphs.
Call $\ul{\cG}$ the subfamily of graphs in $\cG$ that have no multiple edges and $\ul{\cG_1}$, $\ul{\cG_2}$, $\ul{\cG_3}$ the corresponding subfamilies of 
connected, 2-connected, and 3-connected graphs. 

To write down a grammar for $\ul{\cG}$ we just need to trace where the multiple edges might 
appear in the decomposition grammar for $\cG$. Clearly a graph is simple iff all its 2-connected components are simple, so we just need
to look at the last part of the grammar: 2-connected from 3-connected. Among the two 2-connected graphs with less than 3 edges (the family $\epsilon$),
we have to forbid the double-link graph, i.e., we have to take $\epsilon=\ell_2$ instead of $\epsilon=\ell_1+\ell_2$. 
Among the 2-connected graphs with at least 3 edges---those giving rise to a RMT-tree---we have to forbid those where some M-brick has at least 2 components that are edges (indeed all representatives of 
a multiple edge are components of a same M-brick, so we can characterize the absence of multiple edges directly on the RMT-tree).
Accordingly we have to change the specifications of the classes involving the decomposition at an M-brick, i.e., the classes $\cP$, $\cB_M$, and $\cV_M$.
In each case we have to distinguish if there is one edge component or zero edge-component incident to the M-brick, so there are two terms for
the decomposition of each these families.
For the class $\cP$ of parallel networks, we have now 
$\cP=e*\textbf{Set}_{\geq 1}(\cD-\mathcal{P}-e)+ \textbf{Set}_{\geq 2}(\cD-\mathcal{P}-e)$ instead of $\cP=\textbf{Set}_{\geq 2}(\cD-\cP)$.
For the class $\cB_M$, we have now 
$\cB_M=(v^2*e*\Set_{\geq 2}(\cD-\cP-e)+v^2*\Set_{\geq 3}(\cD-\cP-e))/\bullet\leftrightarrows\bullet$
  instead of $\cB_M=(v^2*\Set_{\geq 3}(\cD-\cP))/\bullet\leftrightarrows\bullet$. And for the class $\cV_M$, we have now 
$\cV_M=v * e * \textbf{Set}_{\geq 2}(\cD - \mathcal{P} - e) + v * \textbf{Set}_{\geq 3}(\cD - \mathcal{P} - e)$ instead of $\cV_M=v * \textbf{Set}_{\geq 3}(\cD - \mathcal{P})$. 
 These are indicated in the grammar (Figure~\ref{GrammarFig}).
 
 \begin{figure}
{\footnotesize\begin{tabular}[t]{|l|}
  \hline
  \\[-.2cm]
  \textbf{(1). General from Connected (folklore).} \\[.2cm]

  \hspace{.5cm} $\mathcal{G} = \textbf{Set}(\mathcal{G}_1)$ \\[.2cm]
  
  \hline
    \\[-.2cm]
  \textbf{(2). Connected from 2-Connected (Bergeron, Labelle, Leroux).} \\[.2cm]
  
  \begin{tabular}{r@{\ }c@{\ }l}
  $\cG_1$ &=& $\cC=\cC_v + \cC_B - \cC_{Bv}$ [dissymmetry theorem]\\
  $\cC_v$ &=& $v * \cC'$ \\
  $\cC'$ &=& $\textbf{Set} ( \cGbp\circ_v \cC' )$ \\
  $\cC_B$ &=& $\cG_2\circ_{v} \cC'$ \\
  $\cC_{Bv}$ &=& $(v * \cGbp)\circ_v\cC'$ 
  \end{tabular}
  \\
  \\[-.2cm]
  \hline
  \\[-.2cm]
  \textbf{(3). 2-Connected from 3-Connected.} \\
  \\
  (i) \textbf{Networks} \\[.1cm]
  \begin{tabular}{r@{\ }c@{\ }l}
  $\qquad \cD$ &=& $e + \mathcal{S} + \mathcal{P} + \mathcal{H}$ \\
  $\qquad \mathcal{S}$ &=& $(\cD - \mathcal{S})*v*\cD$ \\[.15cm]
  $\qquad \mathcal{P}$ &=& $\left\{
    \begin{array}{ll}
    \textbf{Set}_{\geq 2}(\cD - \mathcal{P}), & \hbox{[Multi-edges allowed]} \\
    e * \textbf{Set}_{\geq 1}(\cD - \mathcal{P} - e ) + \textbf{Set}_{\geq 2}(\cD - \mathcal{P} - e), & \hbox{[No Multi-edges]}\\
\end{array}%
\right.    $\\[.4cm]
  \qquad$\mathcal{H}$ &=& $\ovr{\cG_3} \circ_{e} \cD$
  \end{tabular}
  \\
  \\
  
  (ii) \textbf{Unrooted 2-Connected} \\[.1cm]
  \hspace{-.7cm}\begin{tabular}{r@{\ }c@{\ }l}
  $\qquad\cG_2$ &=& $\epsilon+\cB$,\ \ \ \ \  $\epsilon=\left\{\begin{array}{ll}\ell_1+\ell_2&\ \hbox{[Multi-edges allowed]}\\ \ell_1&\ \hbox{[No multi-edges]}\end{array}\right.$\\[.1cm]
  $\cB$ &=& $\cB_{R} + \cB_{M} + \cB_{T} - \cB_{R-M} - \cB_{R-T} - \cB_{M-T} - \cB_{T\rightarrow T} + \cB_{T-T}$ \\
  $\qquad\cB_{R}$ &=& $\mathcal{R} \circ_{e} (\cD - \mathcal{S})$ \\
  $\qquad \cB_{M}$ &=& $\left\{%
    \begin{array}{ll}
    \mathcal{M} \circ_{e} (\cD - \mathcal{P})=(v^2*\Set_{\geq 3}(\cD-\cP))/\bullet\leftrightarrows\bullet, & \hbox{[Multi-edges allowed]} \\
    (v^2*e*\Set_{\geq 2}(\cD-\cP-e)+v^2*\Set_{\geq 3}(\cD-\cP-e))/\bullet\leftrightarrows\bullet, & \hbox{[No multi-edges]} \\
    \end{array}%
\right.    $\\[.3cm]
  $\qquad\cB_{T}$ &=& $\cG_3 \circ_{e} \cD$\\
  $\qquad\cB_{R-M}$ &=& $(v^2*\mathcal{S} * \mathcal{P})/\bullet\leftrightarrows\bullet,  \qquad\hbox{[Up to pole exchange, denoted by } /\bullet\leftrightarrows\bullet]$\\
  $\qquad\cB_{R-T}$ &=& $(v^2*\mathcal{S} * \mathcal{H})/\bullet\leftrightarrows\bullet$\\
  $\qquad\cB_{M-T}$ &=& $(v^2*\mathcal{P} * \mathcal{H})/\bullet\leftrightarrows\bullet$\\
  $\qquad\cB_{T\rightarrow T}$ &=& $(v^2*\mathcal{H} * \mathcal{H})/\bullet\leftrightarrows\bullet$\\  $\qquad\cB_{T-T}$ &=& $(v^2*\mathcal{H} * \mathcal{H})/(\bullet\leftrightarrows\bullet , H\leftrightarrows H),  \qquad\hbox{[Up to pole and component exchange]}$
  \end{tabular}
  \\
  \\
  (iii) \textbf{Vertex-pointed 2-Connected} \\[.1cm]
  \hspace{-.7cm}\begin{tabular}{r@{\ }c@{\ }l}  $\qquad \cGbp$ &=& $\epsilon'+\cB'$,\ \ \ \ \  $\epsilon'=\left\{\begin{array}{ll}\ell_1\ \!\!\!'+\ell_2\ \!\!\!'=v*(e+\Set_2(e))&\ \hbox{[Multi-edges allowed]}\\ \ell_1\ \!\!\!'=v*e&\ \hbox{[No multi-edges]}\end{array}\right.$\\[.2cm]
  $\cB'$&=&$\cV=\cV_{R} + \cV_{M} + \cV_{T} - \cV_{R-M} - \cV_{R-T} - \cV_{M-T} - \cV_{T\rightarrow T} + \cV_{T-T}$ \\
  $\qquad\cV_{R}$ &=& $\mathcal{R}'\circ_{e} (\cD - \mathcal{S})$ \\
  $\qquad \cV_{M}$ &=& $\left\{%
    \begin{array}{ll}
    v * \textbf{Set}_{\geq 3}(\cD - \mathcal{P}), & \hbox{[Multi-edges allowed]} \\
    v * e * \textbf{Set}_{\geq 2}(\cD - \mathcal{P} - e) + v * \textbf{Set}_{\geq 3}(\cD - \mathcal{P} - e), & \hbox{[No multi-edges]} \\
    \end{array}%
\right.    $\\
\\[-.3cm]
  $\qquad\cV_{T}$ &=& $ \cGtp \circ_{e} \cD$\\
  $\qquad\cV_{R-M}$ &=& $v*\mathcal{S} * \mathcal{P}$ \\
  $\qquad\cV_{R-T}$ &=& $v*\mathcal{S} * \mathcal{H}$ \\
  $\qquad\cV_{M-T}$ &=& $v*\mathcal{P} * \mathcal{H}$ \\
  $\qquad\cV_{T\rightarrow T}$ &=& $v*\mathcal{H} *\mathcal{H}$\\
    $\qquad\cV_{T-T}$ &=& $(v*\mathcal{H} * \mathcal{H})/H\leftrightarrows H$\\[.1cm]
    \end{tabular}
  \\[.2cm]
  \hline
\end{tabular}}
\caption{The grammar to decompose a graph family $\cG$ stable under Tutte's decomposition. 
The non-basic terminal classes of the grammar are the 3-connected classes $\cG_3$, $\cGtp$, 
and $\vec{\cG_3}$.}
\label{GrammarFig}
\end{figure}

\begin{figure}
{\footnotesize\begin{tabular}[t]{|l|}
  \hline
  \\[-.2cm]
  \textbf{(1). General from Connected (folklore).} \\[.2cm]

  \hspace{.5cm} $ {G}(z,y) = \exp( {G}_1(z,y))$ \\[.2cm]
  
  \hline
    \\[-.2cm]
  \textbf{(2). Connected from 2-Connected (Bergeron, Labelle, Leroux).} \\[.1cm]
  
\hspace{.2cm}{{\large [}The variables $z$ and $x$ are related by $x=zG_1\,\!\!'(z,y)$.{\large ]}}\\[.2cm]
  
  \begin{tabular}{r@{\ }c@{\ }l}
  $ G_1(z,y)$ &=& $ C(z,y)= C_v(z,y) +  C_B(z,y) -  C_{vB}(z,y)$\\
  $ C_v(z,y)$ &=& $zC'(z,y)$ \\
  $ C'(z,y)$ &=& $\exp\left(G_2\,\!\!'(x,y)\right)$ \\
  $C_B(z,y)$ &=& $G_2(x,y)$ \\
  $C_{vB}(z,y)$ &=& $xG_2\,\!\!'(x,y)$ 
  \end{tabular}
  \\
  \\[-.2cm]
  \hline
  \\[-.2cm]
  \textbf{(3). 2-Connected from 3-Connected.} \\[.1cm]
  
  \hspace{.2cm}{{\large [}The variables $y$ and $w$ are related by $w=D(x,y)$.{\large ]}}\\
  \hspace{.2cm}{{\large [}We use the notations $\ds\exp_{\geq k}(t)=\exp(t)-\sum_{i=1}^{k-1}\frac{t^i}{i!}$\ \  and\ \  $\ds{\loga}_{\geq k}(t)=\log(1/(1-t))-\sum_{i=1}^{k-1}\frac{t^i}{i}$.{\large ]}}\\
  
  (i) \textbf{Networks} \\[.1cm]
  \begin{tabular}{r@{\ }c@{\ }l}
  $\quad D(x,y)$ &=& $y +  {S}(x,y) +  {P}(x,y) +  {H}(x,y)$ \\
  $\quad  {S}(x,y)$ &=& $(D(x,y) -  {S}(x,y))x D(x,y)$ \\
  $\quad  {P}(x,y)$ &=& $\left\{
    \begin{array}{ll}
    \exp_{\geq 2}(D(x,y)-P(x,y)) & \hbox{[Multi-edges]} \\
    y\exp_{\geq 1}(D(x,y)-{P}(x,y) - y ) + \exp_{\geq 2}(D(x,y)-{P}(x,y)-y) & \hbox{[No Multi-edges]}\\
\end{array}%
\right.    $\\[.2cm]
  \quad$ {H}(x,y)$ &=& $\ovr{G_3}(x,w)$\\[.3cm]
  \end{tabular}
  \\[-.1cm]
  (ii) \textbf{Unrooted 2-Connected} \\[-.1cm]
  \hspace{-.7cm}\begin{tabular}{r@{\ }c@{\ }l}
  $G_2(x,y)$ &=& $\epsilon(x,y)+B(x,y)$,\ \ \ \ \  $\epsilon(x,y)=\left\{\begin{array}{ll}x^2(y/2+y^2/4)&\ \hbox{[Multi-edges]}\\x^2y/2&\ \hbox{[No multi-edges]}\end{array}\right.$\\[.2cm]
  $B(x,y)$ &=& $B_{R}(x,y)\!+\!B_{M}(x,y)\!+\!B_{T}(x,y)\!-\!B_{R-M}(x,y)\!-\!B_{R-T}(x,y)\!-\!B_{M-T}(x,y)\!-\!B_{T-T}(x,y)$ \\[.1cm]
 $\quad B_{R}(x,y)$ &=& $\ds{\loga}_{\geq 3}(x\cdot(D(x,y)\!-\!S(x,y)))/2$\\
  $\quad B_{M}(x,y)$ &=& $\left\{%
    \begin{array}{ll}
     \ds \!\!\!x^2\exp_{\geq 3}(D(x,y)-{P}(x,y))/2, & \hbox{[Multi-edges]} \\
     \ds \!\!\!x^2\!\left(y\exp_{\geq 2}(D(x,y)\!-\!{P}(x,y)\!-\!y)\!+\!\exp_{\geq 3}(D(x,y)\!-\!{P}(x,y)\!-\!y)\right)\!/2, & \hbox{[No multi-edges]}\!\!\!\!\! \\
    \end{array}%
\right.    $\\[.2cm]
  $\quad B_{T}(x,y)$ &=& $G_3(x,w)$\\
    $\quad \ B_{R\!-\!M}(x,y)$ &=& $x^2{S}(x,y)\cdot{P}(x,y)/2$\\[.0cm]
  $\quad B_{R-T}(x,y)$ &=& $x^2{S}(x,y)\cdot{H}(x,y)/2$\\[.0cm]
  $\quad B_{M\!-\!T}(x,y)$ &=& $x^2{P}(x,y)\cdot{H}(x,y)/2$\\[.0cm]
  $\quad B_{T-T}(x,y)$ &=& $x^2{H}(x,y)^2/4$\\[.2cm]
  \end{tabular}
\\
  (iii) \textbf{Vertex-pointed 2-Connected} \\[-.1cm]
  \hspace{-.7cm}\begin{tabular}{r@{\ }c@{\ }l}  $\qquad G_2\,\!\!'(x,y)$ &=& $\epsilon'(x,y)+B'(x,y)$,\ \ \ \ \  $\epsilon'(x,y)=\left\{\begin{array}{ll}\ \!x(y+y^2/2)&\ \hbox{[Multi-edges]}\\ \ \!xy&\ \hbox{[No multi-edges]}\end{array}\right.$\\[.0cm]
\\[-.4cm]
  $B'(x,y)$&=&$V_{R}(x,y)\!+\!V_{M}(x,y)\!+\!V_{T}(x,y)\!-\!V_{R-M}(x,y)\!-\!V_{R-T}(x,y)\!-\!V_{M-T}(x,y)\!-\!V_{T-T}(x,y)$ \\[.05cm]
  $\quad V_{R}(x,y)$ &=& $\ds x^2(D(x,y)-S(x,y))^2D(x,y)/2$ \\
  $\quad V_{M}(x,y)$ &=& $\left\{%
    \begin{array}{ll}
    x\exp_{\geq 3}(D(x,y)\!-\!{P}(x,y)), & \hbox{[Multi-edges]} \\
    xy\exp_{\geq 2}(D(x,y)\!-\!{P}(x,y)\!-\!y) + x\exp_{\geq 3}(D(x,y)\!-\!{P}(x,y)\!-\!y), & \hbox{[No multi-edges]} \\
    \end{array}%
\right.    $\\
\\[-.4cm]
  $\quad V_{T}(x,y)$ &=& $ G_3\,\!\!'(x,w)$\\[.0cm]
  $\quad \ V_{R\!-\!M}(x,y)$ &=& $x{S}(x,y)\cdot{P}(x,y)$ \\[.0cm]
  $\quad V_{R-T}(x,y)$ &=& $x{S}(x,y)\cdot{H}(x,y)$ \\[.0cm]
  $\quad V_{M\!-\!T}(x,y)$ &=& $x{P}(x,y)\cdot{H}(x,y)$ \\[.0cm]
    $\quad V_{T-T}(x,y)$ &=& $x{H}(x,y)^2/2$\\[.1cm]
    \end{tabular}
  \\
  \hline
\end{tabular}}
\caption{The equation-system to express the series counting a family of graphs in terms
of the series counting the 3-connected subfamilies (the terminal series are $G_3(x,w)$, $G_3\,\!\!'(x,w)$, $\vec{G_3}(x,w)$) is obtained by translating the grammar shown in Figure~\ref{GrammarFig}.}
\label{fig:system}
\end{figure}

\section{Application: counting a graph family reduces to counting the 3-connected subfamily}\label{sec:applica}
At first notice that our grammar makes it possible to get 
an \emph{exact enumeration} procedure for a graph family $\cG$
(stable under Tutte's decomposition) 
from the counting coefficients of the 3-connected
subfamily, whether in the labelled or in the unlabelled setting.
For instance, Figure~\ref{fig:system} shows the translation of the grammar into
an equation system for the corresponding exponential generating functions
(for labelled enumeration). To get the counting coefficients, one can either turn 
this system into recurrences  (see~\cite{bodirsky} for the precise recurrences) 
or simply compute Taylor developments of increasing degree in an inductive
way. For doing all this, a partial decomposition grammar (allowing integrals)
is actually enough, since formal integration just means dividing the $n$th coefficient
by $n$.

The purpose of this section is to stress the relevance of our grammar 
in order to get \emph{analytic expressions} for series counting families of graphs,
thereby making it possible---via the techniques of singularity
analysis---to get precise asymptotic results. (The latter point, singularity analysis,
is not treated in this article, see~\cite{gimeneznoy,GNR07} for a detailed study.) 
Precisely, we show that finding an (implicit) analytic expression for a series counting a family of graphs stable under Tutte's decomposition 
reduces to finding an analytic expression for the series counting the 3-connected subfamily, a task that is easier in many cases, in 
particular for planar graphs, see Section~\ref{sec:planar}. 

Let us mention that Gim\'enez, Noy, and Ru\'e have recently shown  in~\cite{GNR07}
that their method, which involves integrations, also makes it possible to have general equation systems
relating the series counting a family of graphs and the series counting the 3-connected subfamilies.
Hence we do not claim originality here, we only point out that the equation system can be derived automatically from a grammar, without 
using analytic integrations.

\begin{definition}
A multivariate power series $C:=C(\mathbf{z})$  
is said to be \emph{analytically specified} if it is of the following form
(notice that the definition is inductive and relies heavily on substitutions).

\vspace{.1cm}

\noindent \emph{Initial level.} Rational generating functions (in any number of variables),
the $\log$ function, and the $\exp$ function  
are analytically specified.

\vspace{.1cm}

\noindent \emph{Incremental definition}. 
Let $C=F_1$ be a power series that 
is uniquely specified by a system of the form 
\begin{equation}\label{eq:analytic}
\left\{
\begin{array}{rcl}
F_1&=&\mathrm{Expression}_1,\\
\multicolumn{3}{c}{\vdots}\\
F_m&=&\mathrm{Expression}_m,
\end{array}\right.
\end{equation}
each expression being of the form $G(H_1,\ldots,H_k)$,
 where $G$ is a series  \emph{already
known} to be analytically specified, and where each $H_i$ is either equal to one of the
$F_j$'s or is a series already known to be analytically specified.
Then $C$ is declared to be analytically specified as well.
\end{definition}

Analytically specified series are typically amenable to singularity analysis techniques
in order to obtain precise asymptotic informations (enumeration, limit laws of parameters  on a random instance); let us review briefly how the method works.
Given a series $C$ specified by a system of the form~\eqref{eq:analytic},
 one has to consider the dependency graph between the $F_i$'s,
 in particular 
  the strongly connected components of that graph. Notice that 
 some of the $F_i$'s are expressed in~\eqref{eq:analytic} as a substitution in terms of 
 functions already known to be analytically specified; in that case one
 would also include the dependency graph for the specification of this series,
 at a lower level. Therefore one gets a global dependency graph for $C$ 
 that is a hierarchy of dependency graphs (one for each specifying system). 
 
 Then one treats the strongly components of this global dependency graph 
 ``from bottom to top'', tracing the singularities on the way. 
 In general, singularities are either
 due to the Jacobian of the system vanishing (for ``tree-like'' singularities, see 
 the Drmota-Lalley-Woods theorem~\cite[sec VII.6]{FlaSe}) 
 or stem from a singularity of a function in a lower strongly connected component. Recently the techniques of singularity analysis have been 
 successfully applied to the class of planar graphs by Gim\'enez and Noy~\cite{gimeneznoy}, building on earlier results by Bender, Gao, and Wormald~\cite{BeGa}. One important
difficulty in their work was to actually obtain  a system of equations specifying planar graphs without integration operator, since integrations make it difficult to trace the singularities (in particular when several variables are involved).

We prove here that our grammar provides a direct combinatorial way to obtain such a 
system of equations for planar graphs. Thanks to our grammar, 
finding an (implicit) expression for the series counting
a graph family $\cG$ (stable under Tutte's decomposition) reduces to finding an (implicit) expression for the series counting the 3-connected subfamily $\cG_3$.

\begin{theorem}\label{theo:analytic}
Let $\cG$ be a family of graphs satisfying the stability condition (i.e., a graph is in $\cG$ iff its 3-connected components are also in $\cG$).
Call $\cG_3$ the 3-connected subfamily of $\cG$, and
denote respectively by $G(x,y)$ and $G_3(x,y)$ 
the (exponential) counting series for the classes
$\cG$ and $\cG_3$, where $x$ marks the vertices and $y$ marks the edges. Assume that $G_3(x,y)$ is analytically specified. Then $G(x,y)$ is also
analytically specified via the equation-system shown in Figure~\ref{fig:system} and the 
analytic specification for $G_3(x,y)$.
\end{theorem}
\begin{proof}
Assume that $G_3$ is analytically specified. 
At first let us mention the following simple lemma: if a series $C$ is analytically specified
then any of its partial derivatives is also analytically specified; this is due
to the fact that one can ``differentiate'' the specification for $C$. 
(For the sake of illustration, take the univariate series for binary trees:
$B(x)=x+B(x)^2$. Then the derivative $B'(x)$ satisfies the system $\{B'(x)=1+2B'(x)B(x), B(x)=x+B(x)^2\}$.)
Hence,  the series $G_3\,\!\!'$ and $\ovr{G_3}$ are also analytically specified.

Now notice that $G_3$, $G_3\,\!\!'$ and $\ovr{G_3}$ are (together with 
elementary functions such as $\log$) the only terminal
series in the equation-system  shown in Figure~\ref{fig:system}.
In addition, this system clearly satisfies the requirements of an analytic specification. 
One shows successively ``from bottom to top'' that the series for the following 
graph classes are analytically specified: networks, unrooted 2-connected,
vertex-pointed 2-connected, vertex-pointed connected, connected, and general
graphs in the family one considers. 
\end{proof}

\section{The series counting labelled planar graphs can be specified analytically}\label{sec:planar}
This section is dedicated to proving the following result.

\begin{theorem}\label{theo:planar_analytic}
The series counting labelled planar graphs w.r.t. vertices and edges is analytically
specified.
\end{theorem}
Theorem~\ref{theo:planar_analytic}
 has already been proved by Gim\'enez and Noy~\cite{gimeneznoy} (as a crucial
step to enumerate planar graphs asymptotically). Our contribution is a more 
straightforward proof of this result, which uses only combinatorial arguments and elementary algebraic manipulations; 
 the main ingredients are our grammar and a bijective
construction of vertex-pointed maps. In this way we avoid the technical integration steps
addressed in~\cite{gimeneznoy}.

Let us give a more precise outline of our proof. At first we take advantage of the 
grammar; the class of planar graphs
satisfies the stability condition (planarity is preserved by taking the 3-connected 
components), hence proving Theorem~\ref{theo:planar_analytic} reduces (by Theorem~\ref{theo:analytic}) 
to finding an analytic specification
for the series $G_3(x,w)$ 
counting 3-connected planar graphs, which is the task we address 
from now on. 

By a theorem
of Whitney~\cite{Whitney33}, 
3-connected planar graphs have a unique embedding on the sphere (up to 
reflections), hence counting 3-connected planar graphs is equivalent
to counting  3-connected maps (up to reflections). 
To solve the latter problem, we combine several tools. 
Firstly, the Euler relation makes it possible to reduce the enumeration of 3-connected maps
to the enumeration of vertex-pointed 3-connected maps; secondly  
a bijective construction
due to Bouttier et al~\cite{BoDiGu04} allows us to count \emph{unrestricted} 
 vertex-pointed maps; finally a suitable adaptation of our grammar used in a 
 top-to-bottom way yields the enumeration of vertex-pointed 3-connected maps 
 via vertex-pointed 2-connected maps.
 
\subsection{Maps}
A \emph{map} is a connected graph planarly 
embedded on the sphere up to continuous deformation. All maps considered here have at least one edge. As opposed to graphs, loops are allowed on arbitrary maps 
 (as well as multi-edges).
 Thus, under Tutte's definition, a 2-connected map is either the loop-map, 
or is a loopless map with no separating vertex.
Similarly as for graphs, 
a map is rooted if one of its edges is marked and oriented, and a map is vertex-pointed if one of its vertices is marked.
An important remark is that rooting a map is equivalent to marking one of its corners (indeed there is a one-to-one correspondence between the corners and 
the half-edges of a map).
We use the following notations:
\begin{itemize}
\item
$\Mcr(x,s)$ and $\Mcp(x,s)$ are the series counting respectively rooted maps and vertex-pointed maps, where $x$ marks the number of vertices minus one (say, the root vertex for $\Mcr$, the pointed vertex for $\Mcp$), and $s$ marks the number of half-edges. Here vertices are unlabelled and half-edges are labelled (this is enough to avoid any symmetry), hence the series are exponential w.r.t. half-edges and ordinary w.r.t. vertices. 

\item
Similarly, $\Mbr(x,s)$, and $\Mbp(x,s)$ are the series counting rooted 2-connected maps and vertex-pointed 2-connected maps respectively, where $x$ marks the number of (unlabelled) vertices minus one, and $s$ marks the number of (labelled) half-edges.
\item
Finally, $K(x,w)$ ($\Mtr(x,w)$, $\Mtp(x,w)$) are the series of unrooted (rooted, vertex-pointed, resp.) 3-connected maps with at least 4 vertices. 
Contrary to maps and 2-connected maps,
the conventions are exactly the same as for graphs, i.e., $x$ marks labelled vertices (which carry blue labels) and $w$ marks labelled edges (which carry red labels); the marked vertex is unlabelled for maps counted by $\Mtp(x,w)$, and the root-edge and its ends are unlabelled for maps counted by $\Mtr(x,w)$.   
\end{itemize}

\subsection{Reduction to counting vertex-pointed 3-connected maps}\label{sec:reduc_pointed_3conn}
At first we recall Whitney's theorem~\cite{Whi32}: 
every 3-connected planar graph has a unique
embedding on the sphere up to reflection.  Since a 3-connected map with
at least 4 vertices differs from its mirror-image, each planar graph in $\cG_3$ gives rise to exactly two maps, which yields
\begin{equation}
G_3(x,w)=\frac{1}{2}K(x,w).
\end{equation}
Next, the task of finding an expression for $K(x,w)$ is reduced to finding
an expression for the series counting vertex-pointed 3-connected maps.
Let $K^V(x,w)$, $K^E(x,w)$, $K^F(x,w)$ be the series counting respectively
vertex-pointed, edge-pointed, and face-pointed 3-connected maps (here all vertices
and all edges are labelled). The Euler relation (i.e., $1=1/2\cdot(|\mathrm{vertices}|-|\mathrm{edges}|+|\mathrm{faces}|)$) yields
$$
K(x,w)=\frac{1}{2}(K^V(x,w)-K^E(x,w)+K^F(x,w)).
$$
The series $K^E$ corresponds to half the series counting rooted 3-connected planar maps,
which has been obtained by Mullin and Schellenberg~\cite{MuSc68} (we review briefly the computation scheme in Section~\ref{sec:rooted_3_connected}) 
and more recently 
in a bijective way in~\cite{FuPoScL}:
$$
K^E(x,w)=\frac{1}{2}x^2w\ovr{K}(x,w)=\frac1{2}x^2\left(w-\frac{xw^2}{1+xw}-\frac{w^2}{1+w}-\frac{\gamma_1\gamma_2}{xw(1+\gamma_1+\gamma_2)^3}\right),
$$
where $\gamma_1=\gamma_1(x,w)$ and $\gamma_2=\gamma_2(x,w)$ are specified by the system
$$
\gamma_1=xw(1+\gamma_2)^2,\ \ \ \ \gamma_2=w(1+\gamma_1)^2.
$$
Moreover, as the class of 3-connected maps is well known to be stable by duality, we have (using again the Euler  relation)
$$
K^F(x,w)=x^2K^V(1/x,xw).
$$
We conclude that finding an analytic expression for $G_3(x,w)$ reduces to finding
an analytic expression for $K^V(x,w)$, i.e., reduces to finding an analytic expression for $K'(x,w)=K^V(x,w)/x$.

\subsection{Expression for the series of vertex-pointed maps}\label{sec:count_bouttier}
Our goal is now to get an analytic expression for the series $K'(x,w)$ 
counting vertex-pointed 3-connected maps. At first 
we perform this task for the series counting vertex-pointed maps in this section 
(we will go subsequently from vertex-pointed maps to vertex-pointed 3-connected maps), taking advantage of 
 a recent bijective construction.

\subsubsection{Bijection between vertex-pointed maps and mobiles}
The Bouttier, Di Francesco, Guitter bijection, which is described in \cite{BoDiGu04}, relates the problem of counting vertex-pointed
maps to the enumeration of a certain class of trees, called \emph{mobiles}.
We state here the bijection between mobiles and vertex-pointed maps in a slightly
reformulated way.

\begin{definition}
A \emph{mobile} is a plane tree such that:
\begin{itemize}
\item[-] vertices are of two types: white vertices ($\circ$) and black vertices ($\bullet$),
\item[-] edges are of two types: $\bullet-\bullet$ and $\bullet-\circ$,
\item[-] additional ``legs" are pending from black vertices; these legs are not considered as edges,
\item[-] each black vertex has exactly as many pending legs as it has white neighbours.
\end{itemize}
The  $2m$ half-edges (not counting the legs) carry distinct labels in $[1..2m]$.
 \end{definition}
We define $\mathcal{T}[n,2m]$ as the set of mobiles with $n$ (unlabelled) white vertices and $2m$ (labelled)
half-edges, and  
$\mathcal{M}'[n,2m]$ as the set of vertex-pointed maps with $n+1$ (unlabelled) vertices and $2m$ (labelled) half-edges.
The following theorem is a reformulation of~\cite{BoDiGu04}:
\begin{theorem}[Bouttier, Di Francesco and Guitter]
\label{thm:BoDiGu04}
For $n\geq 0$ and $m\geq 1$, there is an explicit  bijection between $\mathcal{M}'[n,2m]$  and $\mathcal{T}[n,2m]$. 
\end{theorem}
The reader familiar with \cite{BoDiGu04} may be surprised by our definition of mobiles: in \cite{BoDiGu04},
the white vertices of the mobiles and the edges of type $\bullet-\bullet$ are labelled by additional positive integers
constrained by certain conditions\footnote{These labels correspond to a distance from a reference vertex, hence they are of a different nature from the labels considered here, which are just used to distinguish the atoms.}. Moreover, the original mobiles of \cite{BoDiGu04} have no legs.
Indeed, the role of the legs is precisely to encode the variations of the constrained labels. 
To wit,  interpreting legs and white neighbours around a black vertex as, respectively, increments and decrements of the  constrained 
labels, one easily reconstructs a valide mobile, in the sense of \cite{BoDiGu04}. We leave the details to the reader.
\begin{figure}
\centerline{\includegraphics[scale=1.0]{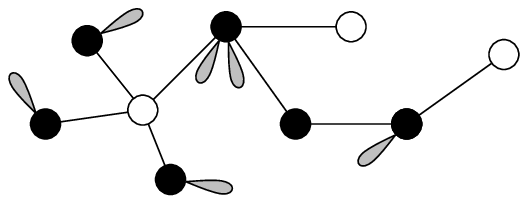}}
\caption{A mobile with $6$ black vertices, $3$ white vertices and $8$ edges.}
\label{fig:mobiles}
\end{figure}

\subsubsection{Families of Motzkin paths}

In order to perform the enumeration of mobiles, we first need to compute some series related to
Motzkin paths.
We consider paths starting at position $0$ and with steps $0$, $+1$, and $-1$.
A bridge is a path ending at position $0$, and an excursion is a bridge that  stays
nonnegative.
We denote respectively by $E$, $B$, and $B^{(+1)}$ the ordinary generating functions of excursions, 
bridges, and paths ending at position $+1$. These series are in two variables, $t$ and $u$,
which count respectively the number of steps $0$ and steps $-1$.
Decomposing these paths at their last passage at $0$, we obtain the following equations:   
\begin{eqnarray*}
E &=& 1 + tE + u E^2, \\
B &=& 1+ (t+2uE)B,  \\
B^{(+1)} &=& EB. 
\end{eqnarray*}
We also need a variant of the series of bridges. Let $\widehat{B}(t,u)$ be the generating function of bridges counted with a weight divided by the total number of steps\footnote{Formally, $\hat{B}$ is obtained from the 
series expansion of $B$ by the substitution $u^kt^\ell \rightarrow \frac{u^kt^\ell}{2k+\ell}$}.
It is easily seen that $\widehat{B}$ is also the series of excursions  counted with a weight
divided by the number of returns to $0$. Hence we have
$$
\widehat{B}(t,u) = \log \frac{1}{1-(t+uE(t,u))} = \log E(t,u).
$$

\subsubsection{Counting families of mobiles}

\begin{figure}
\centerline{\includegraphics[scale=1.0]{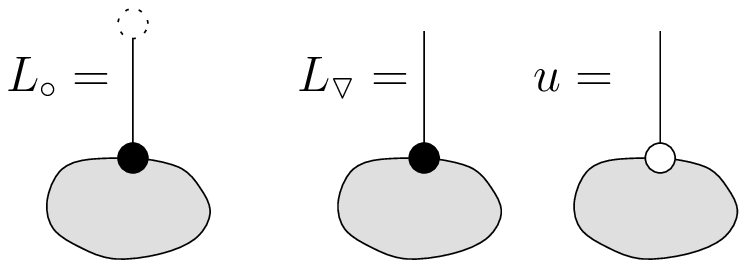}}
\caption{The series $L_\circ$, $L_\triangledown$ and $u$.}
\label{fig:halfmobiles}
\end{figure}

We now compute the generating functions of some classes of mobiles. All the series considered here are exponential in the variable $s$, which counts the number of half-edges  (not counting the legs), and ordinary in the
variables $x$ and $y$, which count respectively the number of white and black vertices.
Following~\cite{BoDiGu04}, we define (see Figure~\ref{fig:halfmobiles}):
\begin{itemize}
\item
$L_{\circ}(x,y,s)$ the  series counting mobiles rooted at an univalent white vertex, which
 is not counted in the series.
\item 
$L_\triangledown(x,y,s)$ the series counting mobiles where an additional edge    
is attached to some black vertex. This edge plays the role of a root  and is counted
in the series (multiplication by $s^2$).
\item
$u(x,y,s)$ the series counting mobiles where an additional edge  
is attached to some white vertex. This edge plays the role of a root; it  
is counted in the series (multiplication by $s^2$) and it has no  vertex at its
other extremity.
\end{itemize}

Observe that one sees exactly a Motzkin bridge when walking around a black vertex (starting at a given position)  
and interpreting white vertices, black vertices, and legs respectively 
as steps $-1$, steps $0$, and steps $+1$. 
Thus a mobile in each type we consider 
is naturally decomposed as a Motzkin bridge 
where each step is substituted by a mobile. 
Precisely we obtain the following equations (see \cite{BoDiGu04}):
\begin{eqnarray}
L_\circ(x,y,s) &=& ys^2 B^{(+1)}\Big(L_\triangledown(x,y,s),u(x,y,s)\Big), \\
L_\triangledown(x,y,s) &=& ys^2 B\Big(L_\triangledown(x,y,s),u(x,y,s)\Big).
\end{eqnarray}
Moreover, one has
\begin{eqnarray}
u(x,y,s) &=& \frac{s^2x}{1-L_\circ(x,y,s)}. 
\end{eqnarray}

All the series involved in the previous equations can be expressed in terms of 
the two simple algebraic series $\beta_1 = \beta_1(x,y,s)$ and $\beta_2 = \beta_2(x,y,s)$
specified by the equation system:
\begin{equation}\label{eq:firstbet}
\left\{
\begin{array}{rcl}
\beta_1 &=& x s^2 +\beta_1^2 + 2 \beta_1 \beta_2, \\[0.1cm]
\beta_2 &=& y s^2 +\beta_2^2 + 2 \beta_1 \beta_2.
\end{array}
\right.
\end{equation}
Then, elementary algebraic manipulations imply
\begin{eqnarray}
usE\Big(L_\triangledown(x,y,s), u(x,y,s)\Big) &=& \beta_1, \\ 
L_\triangledown(x,y,s) &=& \beta_2, \\
E\Big(L_\triangledown(x,y,s),u(x,y,s)\Big)&=& \frac{1}{1-\beta_1-\beta_2}, \\
\frac{1}{1-L_\circ\Big(L_\triangledown(x,y,s),u(x,y,s)\Big)} &=& \frac{\beta_1(1-\beta_1-\beta_2)}{s^2 x}.
\end{eqnarray}

\subsubsection{Counting mobiles}

The enumeration of (unrooted) mobiles is obtained---from the dissymmetry theorem---by 
counting mobiles pointed in several ways. 
We introduce the following series in the three variables $x$, $y$, $s$:
\begin{itemize}
\item[-] $T_\bullet$ ($T_\circ$) is the series of mobiles with a 
distinguished black vertex (white vertex, respectively).
\item[-] $T_{\bullet-\bullet}$ ($T_{\bullet-\circ}$) is the series of mobiles
with a distinguished edge of type $\bullet-\bullet$ ($\bullet-\circ$, respectively). 
\item[-] $T$ is the generating function of all mobiles.
\end{itemize}
The dissymetry theorem implies
\begin{eqnarray}
\label{eq:mobiles}
T(x,y,s) =  T_{\bullet}(x,y,s) + T_{\circ}(x,y,s) - T_{\bullet-\bullet}(x,y,s) - T_{\bullet-\circ}(x,y,s). 
\end{eqnarray}
Moreover, each of the four series above can be expressed in terms of $\beta_1$ and $\beta_2$,
thanks to a decomposition at the pointed vertex or edge:
\begin{eqnarray*}
T_\bullet &=& y \widehat{B}\Big(L_\triangledown(x,y,s), u(x,y,s)\Big) = y \log\frac {1}{1-\beta_1-\beta_2} \\
T_\circ &=& x \log\frac{1}{1-L_\circ\Big(L_\triangledown(x,y,s), u(x,y,s)\Big)} = x\log\frac{\beta_1(1-\beta_1-\beta_2)}{s^2x}\\
T_{\bullet-\bullet} &=& \frac{1}{2s^2} L_{\triangledown}\Big(L_\triangledown(x,y,s), u(x,y,s)\Big)^2 = \frac{1}{2s^2}\beta_2^2 \\
T_{\bullet-\circ} &=& \frac{u}{s^2} L_\circ\Big(L_\triangledown(x,y,s), u(x,y,s)\Big) = \frac{\beta_1\beta_2}{s^2} 
\end{eqnarray*}

In view of Theorem~\ref{thm:BoDiGu04} and Equation~\eqref{eq:mobiles}, we obtain (taking $y=1$):
\begin{lemma}
The series $\Mcp(x,s)$ counting vertex-pointed maps w.r.t. vertices (except one) and half-edges
satisfies the analytic expression 
\begin{equation}\label{eq:Mv}
\Mcp(x,s) =T(x,1,s)= \log \frac{1}{1-\beta_1-\beta_2} + x \log \frac{\beta_1(1-\beta_1-\beta_2)}{s^2x}
-\frac{\beta_2}{2s^2}+\frac{1}{2},
\end{equation}
where $\beta_1$ and $\beta_2$ are specified by~\eqref{eq:firstbet} (with $y=1$).
\end{lemma}

\subsection{Expression\! for\! the\! series\! of\! vertex-pointed 2-connected maps}\label{sec:2conn_maps}
In this section we  take advantage of the block-decomposition to relate the enumeration 
of vertex-pointed maps to the enumeration of vertex-pointed 2-connected maps.
The strategy here differs from the one for graphs (from 2-connected to connected, Section~\ref{grammarCfrom2C}) in two points. Firstly we use the decomposition in a top-to-bottom manner, i.e., we \emph{extract} (instead of expand) the enumeration of vertex-pointed 2-connected maps from the enumeration of vertex-pointed maps obtained in 
Section~\ref{sec:count_bouttier}. Secondly, 
the tree associated here with the decomposition is quite different from the Bv-tree: it takes the embedding into account (the corners play the role of vertices) and it 
is ``completely incident'' to a specific pointed vertex
(similarly as the restricted RMT-tree in Section~\ref{sec:grammarPointed2conn}, which 
applies to 2-connected graphs).

\subsubsection{From rooted maps to rooted 2-connected maps} \label{sec:rooted2conn_fromconn}
Let us first briefly review how a top-to-bottom use of the block-decomposition combined with suitable algebraic manipulations makes it possible to count rooted 2-connected maps from rooted maps, as
demonstrated by Tutte~\cite{Tu63}. Let 
$\ovr{M}(x,s)$ and $\ovr{L}(x,t)$ be the series counting rooted maps and rooted 
2-connected maps w.r.t. vertices (unlabelled) and half-edges (labelled).
Schaeffer~\cite{Sc97} has shown in a bijective way that
\begin{equation}\label{eq:rooted_maps}
\Mcr(x,s)=\frac1{s^4x}\beta_1\beta_2(1-2\beta_1-2\beta_2)-1,
\end{equation}
where $\beta_1=\beta_1(x,s)$ and $\beta_2=\beta_2(x,s)$ are specified by the system
\begin{equation}\label{eq:beta}
\beta_1=xs^2+\beta_1^2+2\beta_1\beta_2,\ \ \ \ \beta_2=s^2+\beta_2^2+2\beta_1\beta_2.
\end{equation}
Given a rooted map $\mu$, its 2-connected core is the maximal 2-connected submap of $\mu$ that contains the root. Clearly, all maps with a given core $\kappa$ 
are obtained by inserting at each corner $c$ of $\kappa$ either a rooted map $\mu_c$ or nothing. Hence $\ovr{M}(x,s)$ and $\ovr{L}(x,t)$ are related by the equation
$$
\Mcr(x,s)=\Mbr(x,s\cdot(1+\Mcr(x,s))).
$$
In other words,
$$
\Mcr(x,s)=\Mbr(x,t),\ \ \ \mathrm{where}\ \ t=s\cdot(1+\Mcr(x,s)).
$$
In view of~\eqref{eq:rooted_maps}, we obtain
$$
t=\frac1{s^3x}\beta_1\beta_2(1-2\beta_1-2\beta_2),
$$ 
hence
$$
xt^2=\frac1{xs^6}\beta_1^2\beta_2^2(1\!-\!2\beta_1\!-\!2\beta_2)^2=\frac{\beta_1^2}{xs^2}\frac{\beta_2^2(1\!-\!2\beta_1\!-\!2\beta_2)^2}{s^4}=\frac{\beta_1}{1\!-\!\beta_1\!-\!2\beta_2}\left(\frac{1\!-\!2\beta_1-2\beta_2}{1-\beta_2-2\beta_1}\right)^2,
$$
and
$$
t^2=\frac1{x^2s^6}\beta_1^2\beta_2^2(1\!-\!2\beta_1\!-\!2\beta_2)^2=\frac{\beta_2^2}{s^2}\frac{\beta_1^2(1\!-\!2\beta_1\!-\!2\beta_2)^2}{x^2s^4}=\frac{\beta_2}{1\!-\!\beta_2\!-\!2\beta_1}\left(\frac{1\!-\!2\beta_1\!-\!2\beta_2}{1\!-\!\beta_1\!-\!2\beta_2}\right)^2.
$$
Equivalently:
$$
xt^2=\eta_1(1-\eta_2)^2,\ \ \ \ t^2=\eta_2(1-\eta_1)^2,
$$
where
\begin{equation}\label{eq:beta_eta}
\eta_1:=\frac{\beta_1}{1-\beta_1-2\beta_2},\ \ \ \ \  \eta_2:=\frac{\beta_2}{1-\beta_2-2\beta_1}.
\end{equation}
The equation-system~\eqref{eq:beta_eta} is readily inverted, giving the following expressions for $\beta_1$ and $\beta_2$ in terms of $\eta_1$ and $\eta_2$:
\begin{equation}\label{eq:eta_beta}
\beta_1:=\frac{\eta_1(1-\eta_2)}{1+\eta_1+\eta_2-3\eta_1\eta_2},\ \ \ \ \  \beta_2:=\frac{\eta_2(1-\eta_1)}{1+\eta_1+\eta_2-3\eta_1\eta_2}.
\end{equation}

Rearranging~\eqref{eq:rooted_maps}, we get the following expression for $\Mcr$ involving only $\beta_1$ and $\beta_2$:
\begin{equation}\label{eq:Mr_no_s}
\Mcr(x,s)= \frac{1-2\beta_1-2\beta_2}{(1-\beta_1-2\beta_2)(1-\beta_2-2\beta_1)}-1,
\end{equation}
hence $\Mbr(x,t)=\Mcr(x,s)$ also satisfies this expression. Replacing $\beta_1$ and $\beta_2$ by their expression in~\eqref{eq:eta_beta}, we obtain the following expression for $\Mbr(x,t)$ in terms of $\eta_1$ and $\eta_2$:
\begin{equation}
\Mbr(x,t)=\eta_1+\eta_2-3\eta_1\eta_2,
\end{equation}
where $\eta_1=\eta_1(x,t)$ and $\eta_2=\eta_2(x,t)$ are specified by the system
\begin{equation}\label{eq:eta}
\eta_1=xt^2(1-\eta_2)^{-2},\ \ \ \ \ \eta_2=t^2(1-\eta_1)^{-2}.
\end{equation}

More generally, if two series $F(x,s)$ and $G(x,t)$ are related by the equation
$$
F(x,s)=G(x,t)\ \ \ \ \mathrm{where}\ \ t=s\cdot(1+\Mcr(x,s)),
$$
and if $F(x,s)$ has an expression in terms of $\beta_1$, $\beta_2$, and $x$, then one gets an expression for $G(x,t)$ in terms of $\eta_1$, $\eta_2$, and $x$ 
 (replacing each occurence of $\beta_1$ or $\beta_2$ according to~\eqref{eq:eta_beta}). This change of variable rule will be useful later on.

\subsubsection{The tree associated to a vertex-pointed map}
\begin{figure}
\centerline{\includegraphics[scale=0.8]{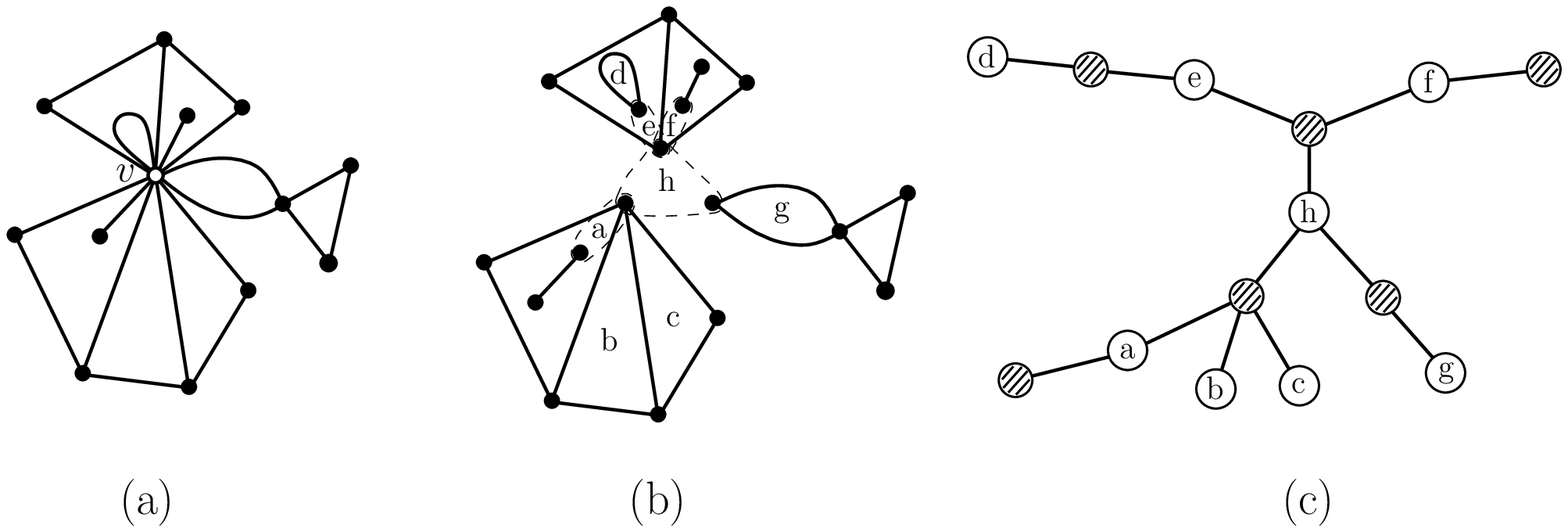}}
\caption{(a) A vertex-pointed map, (b) decomposed according to the blocks and faces  incident to the pointed vertex $v$, (c) the tree resulting from the decomposition (a shaded node
for each block incident to $v$, a letter-node for each face incident to $v$).}
\label{fig:tree_vertex_pointed}
\end{figure}

Let $\mu$ be a vertex-pointed map (as usual, on the sphere), with $v$ the pointed vertex. Let $\mathfrak{B}_v(\mu)$ be the set of (embedded) 2-connected blocks of $\mu$ 
that contain the vertex $v$, and 
let $\mathfrak{F}_v(\mu)$ be the set of faces of $\mu$ that are incident to $v$. 
A face $f\in\mathfrak{F}_v(\mu)$ is said to be incident to a block $b\in\mathfrak{B}_v(\mu)$ if there exists a corner $c$ of $\mu$ incident to $v$, such that the face incident to $c$ is $f$ and 
at least one side of $c$ belongs to $b$.
Call $T$ the bipartite graph with vertex-set  $\mathfrak{F}_v(\mu)\cup\mathfrak{B}_v(\mu)$ and edge-set
corresponding to the incidences faces/blocks. From the Jordan curve theorem and an inductive argument (e.g., on the number of blocks attached at $v$), it is easily checked that $T$ is a tree, see Figure~\ref{fig:tree_vertex_pointed}. 

\subsubsection{From vertex-pointed maps to vertex-pointed 2-connected maps}
We call respectively $\MVB$, $\MVf$, and $\MVBf$ the series counting vertex-pointed maps where, respectively, a  
block-node of $T$ is marked, a face-node of $T$ is marked, and an edge of $T$ is marked. By the dissymmetry theorem (or simply, by the fact that each tree has 
one more vertex than edge), the series $\Mcp$ counting vertex-pointed maps satisfies
\begin{equation}\label{eq:Mcp}
\Mcp=\MVB+\MVf-\MVBf.
\end{equation}
As detailed next, the series $\Mcp$ (obtained in Section~\ref{sec:count_bouttier}),  
and the series $\MVf$ and $\MVBf$ 
have explicit expressions; hence 
Equation~\eqref{eq:Mcp} will yield an expression for $\MVB$, and thus (by the change of variable rule) an
expression for the series $L'$ counting vertex-pointed 2-connected maps. 

Consider a map $\mu$ counted by $\MVf$, i.e., a map with a distinguished incidence face-vertex $(f,v)$.
As $\mu$ might not be  2-connected, the incidence $(f,v)$ can be realised by several corners. 
Let $c_1,\ldots,c_k$ be these corners cyclically ordered in ccw order around $v$, and let $\mu_1,\ldots,\mu_k$
be the components of $\mu$ delimited between successive corners, each of these maps being
naturally rooted at the corner of separation from the other components. 
Clearly the maps $\mu_i$ are exactly constrained to have the incidence root-face/root-vertex realised by a unique corner,
which is their root corner. Such rooted maps are said to have a simple root corner and their series is denoted by $S$.
From the above discussion and  the computation rule associated with the cycle construction (i.e., $\cC=\Cyc(\cA)\Rightarrow C(x)=\log(1/(1-A))$) we get
$$
\MVf=\log\left( \frac1{1-S}\right).
$$
Moreover the series $S$ is easy to obtain from $\Mcr$. Indeed each rooted map is obtained from a rooted map with a single corner where
one may insert another rooted map in the root corner; therefore
$$
\Mcr=S\cdot(1+\Mcr),
$$
i.e., $S=\Mcr/(1+\Mcr)$. Hence 
\begin{equation}\label{eq:Mvf}
\MVf=\log(1+\Mcr).
\end{equation}

Finding an expression for $\MVBf$ in terms of $\Mcr$ is easier. From Figure~\ref{fig:tree_vertex_pointed} 
one easily sees that marking an incidence face-vertex $(f,v)$ + a block
incident to $(f,v)$ is the same as marking a corner of the map. Hence
\begin{equation}\label{eq:Mvbf}
\MVBf=\Mcr.
\end{equation}
Plugging the expression~\eqref{eq:Mr_no_s} 
of $\Mcr$ into~\eqref{eq:Mvf} and~\eqref{eq:Mvbf}, we obtain explicit expressions for $\MVf$ and $\MVBf$ in terms of $\{\beta_1,\beta_2,x\}$.
We can also obtain from~\eqref{eq:Mv} (replacing $s^2$ by $\beta_2(1-\beta_2-2\beta_1)$) an expression for $\Mcp$ in terms of $\{\beta_1,\beta_2,x\}$.

Therefore, we now have expressions for $\Mcp$, $\MVf$ and $\MVBf$ in terms of $\beta_1$, $\beta_2$, and $x$. Plugging these into  $\MVB=\Mcp-\MVf+\MVBf$,
we obtain:
\begin{eqnarray}\label{eq:MV_sans_s}
 \MVB&=&\log\left(\frac{(1-\beta_1-2\beta_2)(1-\beta_2-2\beta_1)}{(1-2\beta_1-2\beta_2)(1-\beta_1-\beta_2)}\right)+x\log\left(\frac{1-\beta_1-\beta_2}{1-\beta_1-2\beta_2}\right)\\
 &&+{\frac {1-3
\,\beta_1-2\,\beta_2}{ 2\left(1- \beta_1-2\, \beta_2\right)  \left( 1- \beta_2-2\, \beta_1 \right) }}-\frac{1}{2},\nonumber
\end{eqnarray}
where $\beta_1$ and $\beta_2$ are specified by~\eqref{eq:beta}.

Clearly each map counted by $\MVB$ arises from a vertex-pointed 2-connected map (the marked 2-connected block) 
where a rooted connected map may be inserted in each corner:
$$
\MVB(x,s)=\Mbp(x,s(1+\Mcr(x,s))).
$$
In other words,
$$
\MVB(x,s)=\Mbp(x,t), \ \ \ \ \mathrm{where}\ \ \ t=s\cdot(1+\Mcr(x,s)).
$$
As described at the end of Section~\ref{sec:rooted2conn_fromconn},  to get an expression for $\Mbp$ in terms of $\eta_1$, $\eta_2$, and $x$, 
we have to take the expression~\eqref{eq:MV_sans_s} of $\MVB(x,s)$, and then replace each occurence of 
$\beta_1$ and $\beta_2$ by their expressions given in~\eqref{eq:eta_beta}. All simplifications done, we obtain
the following expression for the series counting vertex-pointed 2-connected maps: 
\begin{equation}
\Mbp(x,t)=-\log  \left( 1-{ \eta_1}\,{ \eta_2} \right) +x\log  \left( {
\frac {1-{ \eta_1}\,{ \eta_2}}{1-{ \eta_2}}} \right)+\eta_1+\eta_2-3\eta_1\eta_2-\frac{2\eta_1+\eta_2-3\eta_1\eta_2}{2(1-\eta_1)}, 
\end{equation}
where $\eta_1=\eta_1(x,t)$ and $\eta_2=\eta_2(x,t)$ are specified by~\eqref{eq:eta}.

\subsection{Expression\! for\! the\! series\! of\! vertex-pointed 3-connected maps}\label{sec:rooted_3conn}
Our aim in this section is to obtain an analytic specification for $K'$ from the one for $L'$. The strategy is to use the restricted RMT-tree in a top-to-bottom manner
(extract the enumeration of 3-connected  from the enumeration of 2-connected),
while taking the embedding into account. 
\subsubsection{From rooted 2-connected to rooted 3-connected maps}\label{sec:rooted_3_connected}
Similarly as in Section~\ref{sec:2conn_maps}, we first review how rooted 3-connected maps can be extracted from rooted 2-connected maps, following the approach
of Mullin and Schellenberg~\cite{MuSc68}.
We consider here networks that are \emph{embedded} in the plane, with the two poles in the outer face. The series of networks, series-networks,
parallel-networks, and polyhedral networks are still denoted by $\uD=\uD(x,y)$, $\uS=\uS(x,y)$, $\uP=\uP(x,y)$, $\uH=\uH(x,y)$ ($x$ marks vertices not incident to the root, $y$ marks edges different from the root) but in all this section the embedding is taken into account, i.e., we are counting maps. 
Notice that a network is obtained from a rooted 2-connected map with at least 2 edges as follows: distribute distinct blue labels on the vertices not incident to the root, distribute red labels on the edges different from the root, erase the root edge and erase the labels on the half-edges. Therefore
$$
xt^2\uD(x,t^2)=\Mbr(x,t)-(1+x)t^2.
$$ 
In other words (with $y=\sqrt{t}$)
\begin{equation}\label{eq:D_maps}
\uD(x,y)=\frac{\eta_1+\eta_2-3\eta_1\eta_2}{xy}-1/x-1,
\end{equation}
where $\eta_1=\eta_1(x,y)$ and $\eta_2=\eta_2(x,y)$ are specified  by
\begin{equation}
\eta_1=xy(1-\eta_2)^{-2},\ \ \ \ \ \eta_2=y(1-\eta_1)^{-2}.
\end{equation}

Next, Trakhtenbrot's decomposition for rooted 2-connected graphs is readily adapted so as to take the embedding into account; the main difference is that the components of a parallel network 
are ordered (say from left to right if the axis passing by the poles is viewed as vertical). We get the system (equivalent to the one in~\cite{MuSc68}, where everything
is formulated on quadrangulations):

\begin{equation}
\left\{
\begin{array}{rcl}
\uD&=&y+\uS+\uP+\uH,\\
\uS&=&(\uD-\uS)\cdot x\cdot\uD,\\
\uP&=&(\uD-\uP)\cdot\uD,\\
\uH&=&\ovr{K}(x,\uD(x,y)).
\end{array}
\right.
\end{equation}

Hence $$\uH=\uD-\uS-\uP-y=\uD-\frac{x\uD^2}{1+x\uD}-\frac{\uD^2}{1+\uD}-y,$$
which yields
\begin{equation}\label{eq:Kfirst}
\ovr{K}(x,w)=w-\frac{xw^2}{1+xw}-\frac{w^2}{1+w}-y,\ \ \ \mathrm{where}\ \ w=\uD(x,y).
\end{equation}
Replacing $y$ by $\eta_2(1-\eta_1)^2$ and $x$ by $\eta_1(1-\eta_2)^2/(\eta_2(1-\eta_1)^2)$ in the expression~\eqref{eq:D_maps} of $\uD$, we get
$$
w=\frac{\eta_2(1-\eta_1-\eta_2)}{(1-\eta_2)^2}.
$$
and 
$$
xw=\frac{\eta_1(1-\eta_1-\eta_2)}{(1-\eta_1)^2}.
$$
Equivalently, if we introduce
\begin{equation}\label{eq:gamma_eta}
\gamma_1=\frac{\eta_1}{1-\eta_1-\eta_2},\ \ \ \gamma_2=\frac{\eta_2}{1-\eta_1-\eta_2},
\end{equation}
we have
\begin{equation}
xw=\frac{\gamma_1}{(1+\gamma_2)^2},\ \ \ \ \ w=\frac{\gamma_2}{(1+\gamma_1)^2}.
\end{equation}

Now we want to express $\ovr{K}(x,w)$ in terms of $x$, $w$, $\gamma_1$ and $\gamma_2$. In view of Equation~\eqref{eq:Kfirst}, this task is 
clearly equivalent to finding an expression for $y$ in terms of $x$, $w$, $\gamma_1$ and $\gamma_2$. Inverting the system~\eqref{eq:gamma_eta}, we get 
\begin{equation}\label{eq:eta_gamma}
\eta_1=\frac{\gamma_1}{1+\gamma_1+\gamma_2},\ \ \ \ \eta_2=\frac{\gamma_2}{1+\gamma_1+\gamma_2}.
\end{equation}
Substituting these expressions in $y=\eta_2(1-\eta_1)^2$, we obtain
$$
y=\frac{\gamma_2(\gamma_2+1)^2}{(1+\gamma_1+\gamma_2)^3}=\frac{\gamma_1\gamma_2}{xw(1+\gamma_1+\gamma_2)^3}.
$$
This yields the following expression for $\ovr{K}(z,w)$ (see~\cite{MuSc68} for more details):
\begin{equation}
\ovr{K}(z,w)=w-\frac{xw^2}{1+xw}-\frac{w^2}{1+w}-\frac{\gamma_1\gamma_2}{xw(1+\gamma_1+\gamma_2)^3},
\end{equation}
where $\gamma_1=\gamma_1(x,w)$ and $\gamma_2=\gamma_2(x,w)$ are specified by the system
\begin{equation}\label{eq:gamma}
\gamma_1=xw(1+\gamma_2)^2,\ \ \ \ \gamma_2=w(1+\gamma_1)^2.
\end{equation}
More generally, if two series $F(x,y)$ and $G(x,w)$ are related by the equation
$$
F(x,y)=G(x,w),\ \ \mathrm{where}\ w=\uD(x,y),
$$
and if $F(x,y)$ has an expression in terms of $\{x, \eta_1, \eta_2\}$, then an expression for $G(x,w)$ in terms of  $\{x, \gamma_1, \gamma_2\}$ is 
obtained upon replacing each occurence of $\eta_1$ or $\eta_2$ according to~\eqref{eq:eta_gamma}.

\subsubsection{Getting an expression for $K'(x,w)$}\label{sec:exp_pointed_3_conn}
To get an expression for the series $K'(x,w)$ counting vertex-pointed 3-connected maps, the restricted RMT-tree introduced in Section~\ref{sec:root_RMT_tree} 
proves very convenient; it keeps the calculations as simple as possible.
Here we are dealing with maps, so we have to take the embedding into account; this means that the equations relating vertex-pointed 2-connected maps
and vertex-pointed 3-connected maps are the same as for graphs, except that parallel components (M-nodes) are cyclically ordered for maps, due to the embedding.
This yields the following adaptation---taking the embedding into account---of 
the subsystem of equations in Figure~\ref{GrammarFig} related to vertex-pointed 2-connected graphs (here $\uV$
denotes the series counting vertex-pointed 2-connected maps with at least 3 edges, i.e., $\uV(x,y)=L'(x,\sqrt{y})-xy-y/2-xy^2/2$):
$$
\uV=\uV_R+\uV_M+\uV_T-\uV_{R-M}-\uV_{R-T}-\uV_{M-T}-\uV_{T-T},
$$
$$
\uV_R=\frac{1}{2}x^2(\uD-\uS)^2\uD,\ \ \ \ \uV_M=x\ \!\loga_{\geq 3}(D-P),\ \ \ \ \uV_T=K'(x,\uD(x,y)),
$$
where $\loga_{\geq 3}(Z):=\log(1/(1-Z))-Z-Z^2/2$; and 
$$
\uV_{R-M}=x\cdot\uS\cdot\uP,\ \ \ \ \uV_{R-T}=x\cdot\uS\cdot\uH,\ \ \ \ \uV_{M-T}=x\cdot\uP\cdot\uH,\ \ \ \ \uV_{T-T}=\frac1{2}x\cdot\uH^2.
$$
Hence
$$
K'(x,w)=\uV-\uV_R-\uV_M+\uV_{R-M}+\uV_{R-T}+\uV_{M-T}+\uV_{T-T}, \ \ \ \mathrm{where}\  w=D(x,y).
$$
All the series on the right-hand side admit an explicit expression in terms of $\{x,\eta_1,\eta_2\}$. 
Hence $K'(x,w)$ also admits an expression in terms of  $\{x,\eta_1,\eta_2\}$, which is turned into an expression in terms of $\{x,\gamma_1,\gamma_2\}$
when substituting each occurence of $\eta_1$ or $\eta_2$ by the corresponding expression given in~\eqref{eq:eta_gamma}.
(In fact, it is better to allow the presence of $w=D(x,y)$ in all these expressions to have simpler forms, 
similarly as we have done to get an expression 
for $\ovr{K}(x,w)$ in Section~\ref{sec:rooted_3conn}.) All calculations and simplifications done, we obtain the following analytic expression for $K'(x,w)$:

\begin{eqnarray*}
K'(x,w)&=&-\ln  \left( 1-{\frac {{\gamma_1}\,{\gamma_2}}{ \left( 1+{\gamma_1
}+{\gamma_2} \right) ^{2}}} \right) +x\ln  \left( 1+{\frac {{
\gamma_2}\, \left( 1+{\gamma_2} \right) }{ \left( 1+{\gamma_1}+{
\gamma_2} \right)  \left( 1+{\gamma_1} \right) }} \right) \\
 && -{\frac {1\!+\!{\gamma_2}/2\!+\!xw \left( 1\!+\!{\gamma_1} \right) 
 \left( 1\!+\!{\gamma_2} \right)}{1\!+\!{\gamma_1}\!+\!{
\gamma_2}}}-{\frac {3{\gamma_1}\,{\gamma_2}}{ \left( 1\!+\!{
\gamma_1}\!+\!{\gamma_2} \right) ^{2}}}-{\frac {{\gamma_1}\,{
\gamma_2}\, \left( 1\!+\!2\,x\!+\!2\,xw \right) }{2xw \left( 1\!+\!{\gamma_1}\!+\!{
\gamma_2} \right) ^{3}}}\\
&&-\frac1{2x}+1+\frac{w}{2}+xw-\frac{x{w}^{2}}{2}+
{\frac {1}{2\left( 1+xw \right) x}}-x\ln 
 \left( 1+w \right),
\end{eqnarray*}
where $\gamma_1=\gamma_1(x,w)$ and $\gamma_2=\gamma_2(x,w)$ are specified by~\eqref{eq:gamma}.

\subsection{Expression for the series of 3-connected planar graphs}
As we have seen in Section~\ref{sec:reduc_pointed_3conn}, having an expression for the series $K'(x,w)$ (and $\ovr{K}(x,w)$) allows us to obtain 
analytic expressions for the series $G_3(x,w)$, $G_3\,\!\!'(x,w)$, $\ovr{G_3}(x,w)$ counting 3-connected planar graphs, which 
are the terminal series of the equation-system shown in Figure~\ref{fig:system} (assuming we apply this generic equation-system to the family of planar graphs).

\begin{proposition}[expressions for the series counting 3-connected planar graphs]\label{prop:3connected}
The series $\ovr{G_3}(x,w)$, $G_3\,\!\!'(x,w)$, and $G_3(x,w)$ that count respectively   
 rooted, vertex-pointed, and unrooted 3-connected planar graphs w.r.t. vertices 
and edges (recall that one vertex is discarded in $G_3\,\!\!'$, and two vertices and one edge are discarded in $\ovr{G_3}$) satisfy the following analytic expressions
 (the one for  $\ovr{G_3}(x,w)$ was first given in~\cite{MuSc68}, 
 and alternative expressions and computation methods
for $G_3(x,w)$
are given in~\cite{gimeneznoy}):
\begin{equation}
\ovr{G_3}(x,w)=\frac1{2}\left(w-\frac{xw^2}{1+xw}-\frac{w^2}{1+w}-\frac{\gamma_1\gamma_2}{xw(1+\gamma_1+\gamma_2)^3}\right),
\end{equation}
\begin{eqnarray}
G_3\,\!\!'(x,w)&=&-\frac1{2}\ln  \left( 1-{\frac {{\gamma_1}\,{\gamma_2}}{ \left( 1+{\gamma_1
}+{\gamma_2} \right) ^{2}}} \right) +\frac{x}{2}\ln  \left( 1+{\frac {{
\gamma_2}\, \left( 1+{\gamma_2} \right) }{ \left( 1+{\gamma_1}+{
\gamma_2} \right)  \left( 1+{\gamma_1} \right) }} \right)\nonumber\\
 && -{\frac {1+{\gamma_2}/2+xw \left( 1+{\gamma_1} \right) 
 \left( 1+{\gamma_2} \right)}{2(1+{\gamma_1}+{
\gamma_2})}}-{\frac {3{\gamma_1}\,{\gamma_2}}{ 2\left( 1+{
\gamma_1}+{\gamma_2} \right) ^{2}}}-{\frac {{\gamma_1}\,{
\gamma_2}\, \left( 1+2\,x+2\,xw \right) }{4xw \left( 1+{\gamma_1}+{
\gamma_2} \right) ^{3}}} \nonumber \\
&&-\frac1{4x}+\frac1{2}+\frac{w}{4}+\frac{xw}{2}-\frac{x{w}^{2}}{4}+
{\frac {1}{4\left( 1+xw \right) x}}-\frac{x}{2}\ln 
 \left( 1+w \right),
\end{eqnarray}
{\small
\begin{eqnarray}
G_3(x,w)\!\!\!\!\!&=&\!\!\!\!\!\frac{x}{4}\!\left(\!\!-\!2\,\ln\!  \left(1\!-\!{\frac {{\gamma_1}\,{\gamma_2}}{
 \left(\!1\!+\!{\gamma_1}\!+\!{\gamma_2}\! \right) ^{2}}} \right)\!\!+\!x\ln \!
 \left(\!1\!+\!{\frac{{\gamma_2}\, \left(1\!+\! { \gamma_2}\right) }{
 \left(\!1\!+\!{\gamma_1}\!+\!{ \gamma_2} \right)  \left(\!1\!+\!{ \gamma_1}
 \right) }} \right) \!\!+\!\!\frac1{x}\ln \! \left( \!\!1\!+\!{\frac {{
 \gamma_1}\, \left( 1\!+\!{ \gamma_1} \right) }{ \left( 1\!+\!{ \gamma_1}\!+\!{
 \gamma_2} \right)  \left( 1\!+\!{ \gamma_2} \right) }} \right) \right)\nonumber\\
&&\!\!\!\!\!\!\!\!+\frac{x}{4}\left(\!\! -\!\frac1{2}\! -\!\frac{3}{2}\, \frac1{1\!+\! {\gamma_1}\!+\!{  \gamma_2}} \!-\!{\frac {w \left(1\!+\!x \right)  \left( 1\!+\!{  \gamma_1}
 \right)  \left( 1\!+\!{  \gamma_2} \right) }{\!1\!+\!{  \gamma_1}\!+\!{  \gamma_2}}
} \!-\!\,{
\frac {{ 6 \gamma_1}\,{  \gamma_2}}{ \left( \!1\!+\!{  \gamma_1}\!+\!{  \gamma_2}
 \!\right) ^{2}}}\!-\!\frac{3}{2}\,{\frac {{  \gamma_1}\,{  \gamma_2}\, \left( x+1+x
w \right) }{xw \left( 1+{  \gamma_1}+{  \gamma_2}\right) ^{3}}}
   \right)\nonumber\\
&&\!\!\!\!\!\!\!\!+\frac{x}{4}\left(   -\frac1{2}\,x{w}^{2}+xw+w+2   -x\ln  \left( 1+w \right)-{\frac1{x} {\ln  \left( 1+xw \right) }} \right),
\end{eqnarray}
}
where $\gamma_1=\gamma_1(x,w)$ and $\gamma_2=\gamma_2(x,w)$ are specified by 
$$
\gamma_1=xw(1+\gamma_2)^2,\ \ \ \gamma_2=w(1+\gamma_1)^2.
$$
\end{proposition}
\begin{proof}
According to Whitney's theorem, $\ovr{G_3}=\ovr{K}/2$, $G_3\,\!\!'=K'/2$, $G_3=K/2$, where $\ovr{K}$, $K'$ and $K$ are the series counting
rooted, vertex-pointed, and unrooted 3-connected maps. We have already showed how to compute the expressions of $\ovr{K}$ and $K'$ in 
Section~\ref{sec:rooted_3_connected} and Section~\ref{sec:exp_pointed_3_conn}, respectively. Moreover, as discussed 
in Section~\ref{sec:reduc_pointed_3conn},  the series $K$ satisfies
$$
K(x,w)=\frac1{2}x\left( K'(x,w)-\frac{1}{2}xw\ovr{K}(x,w)+K'(1/x,xw)\right).
$$
Notice that the effect on  $(\gamma_1,\gamma_2)$ of replacing $(x,w)$ by $(1/x,xw)$  is simply to swap $\gamma_1$ and $\gamma_2$.
All simplifications done, one obtains the expression of $G_3$ given above.
\end{proof}
\noindent The expressions of Proposition~\ref{prop:3connected}, together with Theorem~\ref{theo:analytic}, finish the proof of Theorem~\ref{theo:planar_analytic}.


\subsection*{Acknowledgements.} The authors thank  Pierre Leroux and Konstantinos Panagiotou for very interesting discussions.

\bibliographystyle{plain}
\bibliography{mabiblio}

\begin{thebibliography}{10}

\bibitem{BeGa}
E.~Bender, Z.~Gao, and N.~Wormald.
\newblock The number of labeled 2-connected planar graphs.
\newblock {\em Electron. J. Combin.}, 9:1--13, 2002.

\bibitem{BeLaLe}
F.~Bergeron, G.~Labelle, and P.~Leroux.
\newblock {\em Combinatorial Species and Tree-like Structures}.
\newblock Cambridge University Press, 1997.

\bibitem{BoFuKaVi07b}
M.~Bodirsky, \'E. Fusy, M.~Kang, and S.~Vigerske.
\newblock Enumeration and asymptotic properties of unlabeled outerplanar
  graphs.
\newblock {\em Elect. J. Comb.}, 14, R66, Sep 2007.

\bibitem{BoGiKaNo07}
M.~Bodirsky, O.~Gim\'{e}nez, M.~Kang, and M.~Noy.
\newblock Enumeration and limit laws for series-parallel graphs.
\newblock {\em Eur. J. Comb.}, 28(8):2091--2105, 2007.

\bibitem{bodirsky}
M.~Bodirsky, C.~Groepl, and M.~Kang.
\newblock Generating labeled planar graphs uniformly at random.
\newblock {\em Theoretical Computer Science}, 379:377--386, 2007.

\bibitem{BKLM07}
M.~Bodirsky, M.~Kang, M.~L\"{o}ffler, and C.~McDiarmid.
\newblock Random cubic planar graphs.
\newblock {\em Random Structures and Algorithms}, 30:78--94, 2007.

\bibitem{BoDiGu04}
J.~Bouttier, P.~Di Francesco, and E.~Guitter.
\newblock Planar maps as labeled mobiles.
\newblock {\em Electr. J. Comb.}, 11(1), 2004.

\bibitem{Diestel}
R.~Diestel.
\newblock {\em Graph Theory}.
\newblock Springer-Verlag, 3rd edition edition, 2005.

\bibitem{FlaSe}
P.~Flajolet and R.~Sedgewick.
\newblock Analytic combinatorics.
\newblock Preliminary version available at
  \texttt{http://algo.inria.fr/flajolet/Publications}.

\bibitem{Fu05a}
{\'E}.~Fusy.
\newblock Quadratic exact size and linear approximate size random generation of
  planar graphs.
\newblock {\em Discrete Mathematics and Theoretical Computer Science},
  AD:125--138, 2005.

\bibitem{FuPoScL}
\'E. Fusy, D.~Poulalhon, and G.~Schaeffer.
\newblock Dissections, orientations, and trees, with applications to optimal
  mesh encoding and to random sampling.
\newblock {\em Transactions on Algorithms}, 4(2):Art. 19, April 2008.

\bibitem{GLL07}
A.~Gagarin, G.~Labelle, and P.~Leroux.
\newblock Counting unlabelled toroidal graphs with no $k_{3,3}$-subdivisions.
\newblock {\em Advances in Applied Mathematics}, 39:51--74, 2007.

\bibitem{GLLW08}
A.~Gagarin, G.~Labelle, P.~Leroux, and T.~Walsh.
\newblock Two-connected graphs with prescribed three-connected components.
\newblock arXiv:0712.1869, 2007.

\bibitem{GNR07}
O.~Gim\'enez, M.~Noy, and J.~Ru\'{e}.
\newblock Graph classes with given 3-connected components: asymptotic counting
  and critical phenomena.
\newblock {\em Electronic Notes in Discrete Mathematics}, 29:521--529, 2007.

\bibitem{gimeneznoy}
Omer Gim\'{e}nez and Marc Noy.
\newblock Asymptotic enumeration and limit laws of planar graphs.
\newblock arXiv:math/0501269, 2005.

\bibitem{Ha}
F.~Harary and E.~Palmer.
\newblock {\em Graphical Enumeration}.
\newblock Academic Press, New York, 1973.

\bibitem{HoTa73}
J.~E. Hopcroft and R.~E. Tarjan.
\newblock Dividing a graph into triconnected components.
\newblock {\em SIAM J. Comput.}, 2(3):135--158, 1973.

\bibitem{Le88}
P.~Leroux.
\newblock Methoden der {A}nzahlbestimmung fur einige {K}lassen von {G}raphen.
\newblock {\em Bayreuther Mathematische Schriften}, 26:1--36, 1988.

\bibitem{MSW05}
C.~McDiarmid, A.~Steger, and D.~Welsh.
\newblock Random planar graphs.
\newblock {\em J. Comb. Theory, Ser. B}, 93:187--205, 2005.

\bibitem{Mohar}
B.~Mohar and C.~Thomassen.
\newblock {\em Graphs on Surfaces}.
\newblock Johns Hopkins University Press, 2001.

\bibitem{MuSc68}
R.C. Mullin and P.J. Schellenberg.
\newblock The enumeration of c-nets via quadrangulations.
\newblock {\em J. Combin. Theory}, 4:259--276, 1968.

\bibitem{Ot48}
R.~Otter.
\newblock The number of trees.
\newblock {\em Annals of Mathematics}, 49(3):583--599, 1948.

\bibitem{Ro70}
R.W. Robinson.
\newblock Enumeration of non-separable graphs.
\newblock {\em J. Comb. Theory}, 9:327--356, 1970.

\bibitem{Sc97}
G.~Schaeffer.
\newblock Bijective census and random generation of {E}ulerian planar maps with
  prescribed vertex degrees.
\newblock {\em Electron. J. Combin.}, 4(1):\string# 20, 14 pp., 1997.

\bibitem{Sh07}
B.~Shoilekova.
\newblock Unlabelled enumeration of cacti graphs.
\newblock Manuscript, 2007.

\bibitem{trak}
B.~A. Trakhtenbrot.
\newblock Towards a theory of non--repeating contact schemes (russian).
\newblock In {\em Trudi Mat. Inst. Akad. Nauk SSSR 51}, pages 226--269, 1958.

\bibitem{Tu63}
W.~T. Tutte.
\newblock A census of planar maps.
\newblock {\em Canad. J. Math.}, 15:249--271, 1963.

\bibitem{Tutte}
W.T. Tutte.
\newblock {\em Connectivity in graphs}.
\newblock Oxford U.P, 1966.

\bibitem{Wa82a}
T.~R.~S. Walsh.
\newblock Counting labelled three-connected and homeomorphically irreducible
  two-connected graphs.
\newblock {\em J. Combin. Theory}, 32(B):1--11, 1982.

\bibitem{Wa82b}
T.~R.~S. Walsh.
\newblock Counting unlabelled three-connected and homeomorphically irreducible
  two-connected graphs.
\newblock {\em J. Combin. Theory}, 32(B):12--32, 1982.

\bibitem{Whi32}
H.~Whitney.
\newblock Congruent graphs and the connectivity of graphs.
\newblock {\em Amer. J. Math.}, 54:150--168, 1932.

\bibitem{Whitney33}
H.~Whitney.
\newblock 2-isomorphic graphs.
\newblock {\em Amer. J. Math.}, 54:245--254, 1933.

\end{thebibliography}
\end{document}